\RequirePackage{etoolbox}
\csdef{input@path}{%
 {sty/}
 {img/}
}%
\csgdef{bibdir}{bib/}

\documentclass[ba,preprint]{imsart}
\pubyear{0000}
\volume{00}
\issue{0}
\doi{0000}
\firstpage{1}
\lastpage{1}

\usepackage{amsthm}
\usepackage{amsmath}
\usepackage{natbib}
\usepackage[colorlinks,citecolor=blue,urlcolor=blue,filecolor=blue,backref=page]{hyperref}
\usepackage{graphicx}

\usepackage{amssymb,bm,bbm}
\usepackage{tikz}
\usepackage{graphicx}
\usepackage{comment}
\usepackage{mathtools}
\usepackage{rotating}
\usepackage{chngcntr}
\usepackage{apptools}
\usepackage[titletoc,title]{appendix}

\usepackage{changebar}

\startlocaldefs
\DeclareMathOperator*{\argmin}{arg\,min}
\DeclareMathOperator*{\argmax}{arg\,max}
\DeclareMathOperator*{\arginf}{arg\,inf}

\newcommand{\transpose}{^{\mathrm{T}}}
\newcommand{\lsigma}{{\underline{\sigma}}}
\newcommand{\usigma}{{\overline{\sigma}}}
\newcommand{\llh}{{\underline{h}}}
\newcommand{\uh}{{\overline{h}}}
\newcommand{\leps}{{\underline{\epsilon}}}

\newcommand{\bdot}{\boldsymbol{\cdot}}
\newcommand{\calA}{{\mathcal{A}}}

\newcommand{\calD}{{\mathcal{D}}}

\newcommand{\calF}{{\mathcal{F}}}
\newcommand{\calG}{{\mathcal{G}}}

\newcommand{\calM}{{\mathcal{M}}}
\newcommand{\calN}{{\mathcal{N}}}

\newcommand{\calP}{{\mathcal{P}}}

\newcommand{\calX}{{\mathcal{X}}}

\newcommand{\calZ}{{\mathcal{Z}}}

\newcommand{\bx}{{\mathbf{x}}}

\newcommand{\bk}{{\mathbf{k}}}

\newcommand{\bl}{{\mathbf{l}}}

\newcommand{\bs}{{\mathbf{s}}}

\newcommand{\bz}{{\mathbf{z}}}

\newcommand{\bbeta}{{\bm{\beta}}}
\newcommand{\bDelta}{{\bm{\Delta}}}

\newcommand{\eye}{{\mathbf{I}}}

\newcommand{\bmu}{{\bm{\mu}}}

\newcommand{\zero}{{\bm{0}}}

\newcommand{\eps}{\epsilon}

\def\hei{\color{black}}

\newcommand{\xx}{\hei \rm}

\newtheorem{proposition}{Proposition}
\newtheorem{theorem}{Theorem}

\newtheorem{lemma}{Lemma}

\theoremstyle{definition}

\newtheorem{remark}{Remark}

\theoremstyle{remark}

\endlocaldefs

\begin{document}


\begin{frontmatter}
\title{Adaptive Bayesian nonparametric regression using a kernel mixture of polynomials with application to the partial linear model}

\runtitle{Kernel Mixture of Polynomials}

\begin{aug}
\author{\fnms{Fangzheng} \snm{Xie}\thanksref{addr1} \ead[label=e1]{fxie5@jhu.edu}}
\and
\author{\fnms{Yanxun} \snm{Xu}\thanksref{addr1,addr2}\ead[label=e2]{yanxun.xu@jhu.edu}}

\runauthor{F. Xie and Y. Xu}

\address[addr1]{Department of Applied Mathematics and Statistics, Johns Hopkins University}
\address[addr2]{Correspondence should be addressed to 
				\printead{e2}}


\end{aug}

\begin{abstract}
We propose a kernel mixture of polynomials prior for Bayesian nonparametric regression. The regression function is modeled by local averages of polynomials with kernel mixture weights. We obtain the minimax-optimal rate of contraction of the full posterior distribution up to a logarithmic factor that adapts to the smoothness level of the true function by estimating metric entropies of certain function classes. We also provide a frequentist sieve maximum likelihood estimator with a near-optimal convergence rate. We further investigate the application of the kernel mixture of polynomials to the partial linear model and obtain both the near-optimal rate of contraction for the nonparametric component and the Bernstein-von Mises limit (i.e., asymptotic normality) of the parametric component. 
The proposed method is illustrated with numerical examples and shows superior performance in terms of computational efficiency, accuracy, and uncertainty quantification compared to the local polynomial regression, DiceKriging, and the robust Gaussian stochastic process.
\end{abstract}


\end{frontmatter}

\section{Introduction}
\label{sec:introduction}
  The standard nonparametric regression model is of the form $y_i=f(\bx_i)+e_i$, where $y_i$'s are observations at given design points, $\bx_i$'s are in the design space $\mathcal{X}\subset\mathbb{R}^p$, and $e_i$'s are independently $\mathrm{N}(0,\sigma^2)$ distributed noises, $i=1,\ldots,n$. The inference task is to  estimate the unknown function 
  $f: \calX\to\mathbb{R}$. Nonparametric regression methods have been widely used in a variety of applications, such as pattern recognition \citep{gyorfi2006distribution,devroye2013probabilistic}, image processing and reconstruction \citep{takeda2007kernel}, electronic healthcare records \citep{xuyanbo2016bayesian}, and semiparametric econometrics \citep{robinson1988root,klein1993efficient}. 
  
  Frequentist methods for nonparametric regression typically compute a fixed estimated function through the given data $(\bx_i,y_i)_{i=1}^n$. 
  In contrast, Bayesian nonparametric techniques first impose a carefully-selected prior distribution on the unknown function $f$ and then find the posterior distribution of $f$ given the observed data $(\bx_i,y_i)_{i=1}^n$, providing a natural way for uncertainty quantification through the full posterior distributions
  instead of a point estimate given by frequentist approaches. 
  One of the most popular Bayesian nonparametric regression methods is the Gaussian process \citep{rasmussen2006gaussian} due to its  tractability. 
  Nevertheless the computational burden of the Gaussian process in likelihood function evaluation resulting from the inversion of the covariance matrix prevents its scalability to big data. 

  In this paper, we propose a novel prior model for nonparametric regression, called \emph{the kernel mixture of polynomials}, that features attractive theoretical properties, efficient computation, and flexibility for extension. Theoretically, we show that by using the kernel mixture of polynomials for Bayesian nonparametric regression, the rate of contraction with respect to the $L_2$-topology is minimax-optimal \citep{stone1982optimal,gyorfi2006distribution} (up to a logarithmic factor) 
    and is adaptive to the smoothness level in the sense that the prior does not depend on the smoothness level of the true function.  
    It is worth mentioning that most papers concerning posterior convergence for Bayesian nonparametric regression only discuss the rate of contraction with respect to the weaker empirical $L_2$-norm \citep{van2008rates,de2010adaptive,van2009adaptive,bhattacharya2014anisotropic}, {i.e.}, the convergence of the function at the given design points. There is little discussion about the rate of contraction with respect to the exact $L_2$-norm for general Bayesian nonparametric regression methods.
    \cite{vaart2011information}, \cite{yoo2016supremum}, and \cite{yang2017frequentist} address this issue only in the context of Gaussian process regression. 
    In particular, the rate of contraction for Gaussian processes with respect to the exact $L_2$-norm requires prior knowledge of the smoothness level of the  true regression function. 
    We also obtain a sieve maximum likelihood estimator with the near-optimal convergence rate as a frequentist point estimator, which could potentially be useful for designing scalable optimization algorithms. 

    From the computational perspective, the proposed kernel mixture of polynomials model 
    avoids the cumbersome $O(n^3)$ inversion of large covariance matrices, greatly facilitating computational efficiency for posterior inference compared to Gaussian process priors, while maintaining the same level of accuracy (see Section \ref{sec:comparison_with_rescaled_gaussian_processes} for detailed comparisons). Such a nice computational advantage makes it attractive to the big-data regime. The code for implementation is publicly available at \url{https://github.com/fangzhengxie/Kernel-Mixture-of-Polynomials.git}. 

In addition, the kernel mixture of polynomials is flexible for extension due to its nonparametric nature. As a specific example, we study the application of this prior model to the partial linear model. The partial linear model is a classical semiparametric regression model of the form $y_i=\bz_i\transpose{}\bbeta+\eta(\bx_i)+e_i$, where $\bz_i,\bx_i$'s are design points, $\bbeta$ is the linear coefficient, $\eta$ is some unknown function, and $e_i$'s are independent $\mathrm{N}(0,1)$ noises, $i=1,\ldots,n$. The literatures of partial linear models from both the frequentist perspective \citep{engle1986semiparametric,chen1988convergence,speckman1988kernel,hastie1990generalized,fan1999root} and Bayesian approaches \citep{lenk1999bayesian,bickel2012semiparametric,tang2015bayesian,yang2015semiparametric} are rich. However, there is little discussion regarding the theoretical properties of the Bayesian partial linear model. 
  To the best of our knowledge, only \cite{bickel2012semiparametric} and \cite{yang2015semiparametric}  discuss the asymptotic behavior of the marginal posterior distribution of $\bbeta$ with Gaussian process priors on $\eta$. We impose the kernel mixture of polynomials prior on $\eta$ and obtain both a near-optimal rate of contraction for $\eta$ and the Bernstein-von Mises limit (i.e., asymptotic normality) of the marginal posterior of $\bbeta$. 

The layout of this paper is as follows. Section \ref{sec:preliminaries} presents the setup for the kernel mixture of polynomials prior for nonparametric regression. Section \ref{sec:theoretical_properties_of_the_bayesian_kernel_regression} elaborates on the convergence properties of the kernel mixture of polynomials for nonparametric regression. 
Section \ref{sec:bayesian_kernel_smoothing_in_partial_linear_models} presents the application of the kernel mixture of polynomials to the partial linear model. 
Section \ref{sec:numerical_studies} illustrates the proposed methodology using numerical examples. We conclude the paper with several discussions in Section \ref{sec:discussion}. 

\section{Preliminaries} 
\label{sec:preliminaries}

\subsection{Setup} 
\label{sub:bayesian_kernel_smoother}
   Recall that the Gaussian nonparametric regression is of the form
   \begin{align}
   y_i&=f(\bx_i)+e_i,\nonumber\\
   \bx_i&\sim p_\bx\text{ independently},\nonumber\\
   e_i&\sim\mathrm{N}(0,\sigma^2)\text{ independently},\quad i = 1,\ldots,n,\nonumber
   \end{align}
   where $p_\bx$ is the marginal density of the design points $(\bx_i)_{i = 1}^n$ supported on the $p$-dimensional unit hypercube $\calX = [0, 1]^p$. We assume that $p_\bx$ is bounded away from $0$ and $\infty$, and is known and fixed.
   The true but unknown regression function $f_0$ is assumed to be in the $\alpha$-H\"older function class $\mathfrak{C}^{\alpha,B}(\calX)$ with envelope $B$, the class of functions $f$ that are $\lceil\alpha -1\rceil$-times continuously differentiable with 
  \begin{align}
  \max_{|\bs|\leq\lceil\alpha-1\rceil}\|D^\bs f\|_\infty+\max_{|\bs|=\lceil\alpha-1\rceil}\sup_{\bx_1\neq \bx_2}\frac{|D^\bs f(\bx_1)-D^\bs f(\bx_2)|}{\|\bx_1-\bx_2\|^{\alpha-|\bs|}}\leq B\nonumber
  \end{align}
  for all $\bx_1,\bx_2\in\calX$, where $\lceil \alpha-1\rceil$ denotes the minimum integer no less than $\alpha-1$, 
  $D^\bs=\partial^{|\bs|}/\partial x_1^{s_1}\ldots\partial x_p^{s_p}$ is the mixed partial derivative operator, and $|\bs|=\sum_{j=1}^ps_j$. When $|\bs|=0$, by convention we define $D^\zero f(\bx)=f(\bx)$.
  We use $\mathbb{P}_0$ and $\mathbb{E}_0$ to denote the probability and expected value under $p_0$, respectively. For readers' convenience, descriptions of additional notations are provided in Appendix. 

  The goal is to estimate the unknown function $f_0$. Leaving the Bayesian framework for a moment, let us consider the frequentist Nadaraya-Watson estimator \citep{nadaraya1964estimating,watson1964smooth} of the form
  \begin{align}\label{eqn:NW_estimator}
  \widehat{f}(\bx)=\sum_{i = 1}^n\left[\frac{\varphi_h(\bx-\bx_i)}{\sum_{i = 1}^n\varphi_h(\bx-\bx_i)}\right]y_i,
  \end{align}
  where $\varphi_h:\mathbb{R}^p\to [0,+\infty)$ is the kernel function parametrized by the bandwidth parameter $h\in(0,+\infty)$ and is assumed to decrease when $\|\bx\|$ increases. It is a local averaging estimator \citep{gyorfi2006distribution}, since the summand ${\varphi_h(\bx-\bx_i)}/{\sum_i\varphi_h(\bx-\bx_i)}$ can be treated as the weight received from $y_i$. 
  As such a simple local averaging estimator does not yield an optimal rate of convergence when the true regression function is $\alpha$-H\"older for $\alpha\geq2$ \citep{devroye2013probabilistic}, \cite{fan1996local} considers a more general local polynomial regression to capture higher-order curvature information of the unknown regression function and gain an optimal rate of convergence. 
  Inspired by these two classical approaches for nonparametric regression, we develop the kernel mixture of polynomials model.

  Firstly, for a given integer $K$, we partition the design space $\calX$ into $K^p$ disjoint hypercubes: 
     $\calX=\bigcup_{\bk\in[K]^p}\calX_K({\bk})$, where 
    $\calX_{K}({\bk})=\prod_{j=1}^p\left({(k_j-1)}/{K}, {k_j}/{K}\right]$, ${\bmu}_{\bk}^\star=\left[{(2k_1-1)}/{(2K)},\cdots,{(2k_p-1)}/{(2K)}\right]\transpose{}$,
   ${\bk}=[k_1,\ldots,k_p]\transpose{}\in[K]^p$, and $\bmu_\bk^\star$ is the center of $\calX_K({\bk})$. The idea of the partition is to breakdown the problem of estimating the regression function over the entire domain $\calX$ into sub-problems of estimating the regression function within each block. 

\begin{figure}[!ht] 
    \centerline{\includegraphics[width=.7\textwidth]{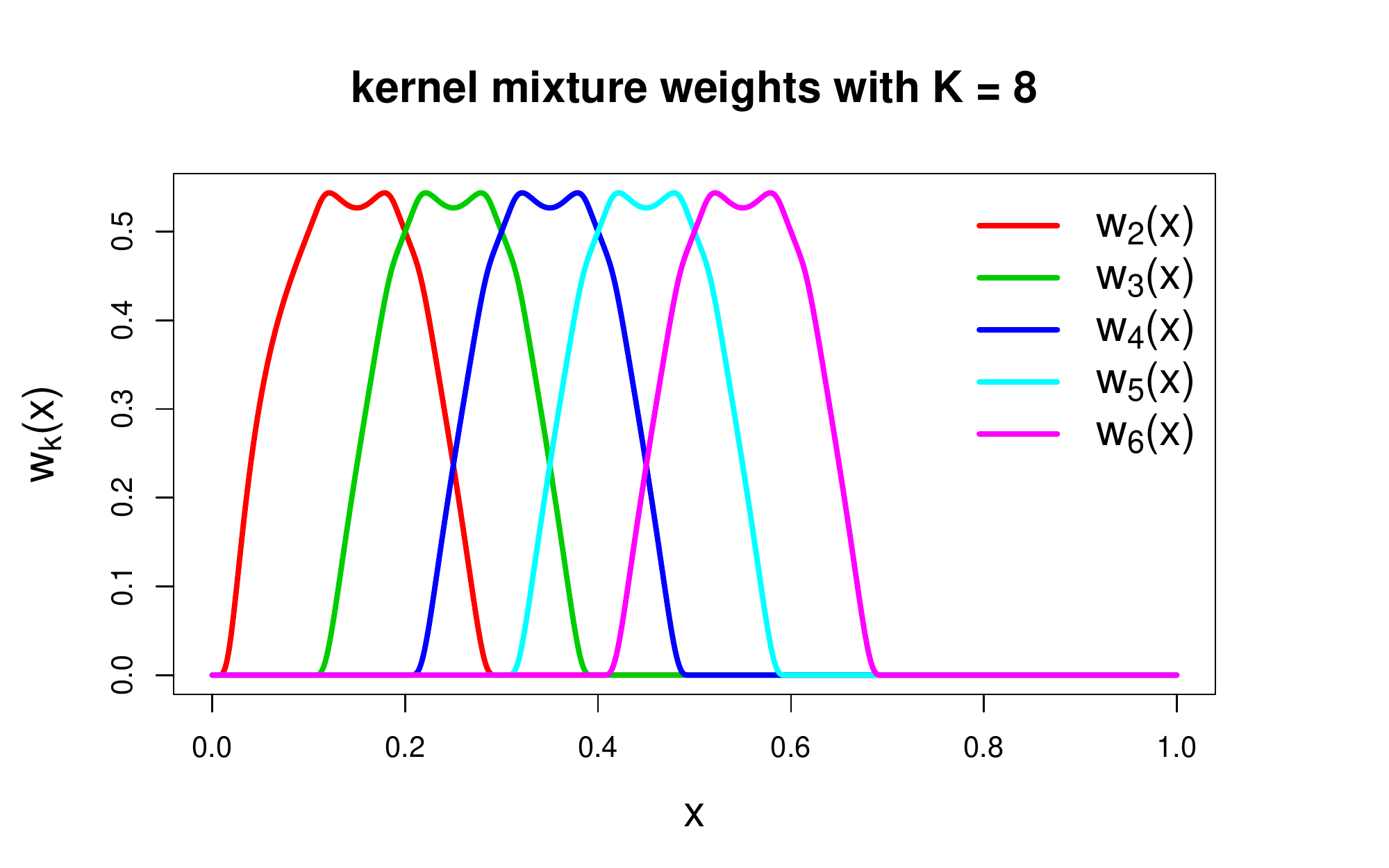}}
    \caption{Examples of $w_k(x)$ for $k = 2, 3, 4, 5, 6$ with $K = 8$ and equidistant $\mu_k = (2k - 1)/(2K)$, $k\in[K]$. }
    \label{Fig:kernel_mixture_weights}
    \end{figure}

    Next, we introduce the notion of (boxed) kernel functions in order to detect the local behavior of the underlying regression function $f_0$  around the block $\calX_K(\bk)$. Formally, a continuous function $\varphi:\mathbb{R}^p\to[0,1]$ is a boxed kernel function if it is supported on $\{{\bx}:\|{\bx}\|_\infty\leq 1\}$, does not increase when $\|{\bx}\|_\infty$ increases, and 
        $\varphi({\bx})\leq\mathbf{1}(\|{\bx}\|_\infty\leq 1)$. 
        We consider the
    univariate bump kernel $\varphi(x)=\exp[-(1-x^2)^{-1}]\mathbf{1}(|{x}|<1)$ in this paper. Other examples of boxed kernels include the triangle kernel, the Epanechnikov kernel, etc. For convenience we denote $\varphi_h({\bx}) = \varphi({\bx}/h)$, where $h>0$ is the {bandwidth} parameter. 
    For each $K\in\mathbb{N}_+$ and each ${\bk}\in[K]^p$, define the kernel mixture weight as
    \begin{align}\label{eqn:kernel_mixture_weights}
    w_{\bk}({\bx})=\frac{\varphi_h({\bx}-{\bmu}_{\bk})}{\sum_{{\bl}\in[K]^p}\varphi_h({\bx}-{\bmu}_{\bl})},
    \end{align} 
    where ${\bmu}_{\bk}\in\calX_K({\bk})$, and $h>0$. 
    The kernel mixture weight $ w_{\bk}({\bx})$ is motivated by the form of Nadaraya-Watson estimator \eqref{eqn:NW_estimator} and is designed to average the signal from $f_0$ locally around $\bmu_\bk$. 
    It can also be shown that the denominator $D({\bx}):=\sum_{{\bl}}\varphi_h({\bx}-{\bmu}_{\bl})$ is strictly non-zero by considering the structure of the disjoint blocks $\calX_K(\bk)$'s. See Figure \ref{Fig:kernel_mixture_weights} for some examples of the kernel mixture weights $w_k(x)$ in the univariate case.  

Thirdly, in order to enhance the estimation accuracy, we need additional tools to capture the higher-order curvature information (\emph{e.g.}, gradient, Hessian, etc.) of the underlying regression function surface, motivated by the local polynomial regression. Given $K$ and  the kernel mixture weights \eqref{eqn:kernel_mixture_weights}, we define the {kernel mixture of polynomials system} to be the set of functions of the form
    \begin{align}
    \label{eqn:KEMPOL}
    \psi_{\bk\bs}(\bx)=w_\bk(\bx)(\bx-\bmu_\bk^\star)^\bs,\quad \bk\in[K]^p,\quad \bs\in\{\bs\in[m]^p:|\bs|=0,1,\ldots,m\},
    \end{align}
    where $(\bx-\bmu_\bk^\star)^\bs$ denotes the monomial $\prod_{j=1}^p[x_j-(2k_j-1)/2K]^{s_j}$, $\bx = [x_1,\cdots,x_p]\transpose{}$, and $\bs = [s_1,\cdots,s_p]\transpose{}$. 
    In this work we require that the degree $m$ of the kernel mixture of polynomials system to be no less than $\alpha$, {i.e.}, $m\geq\alpha$. 
    A similar assumption was also adopted in \cite{de2012adaptive} in the context of conditional Gaussian tensor-product spline models, where the degree of polynomials in the B-splines was required to be no less than the smoothness level of the underlying true function. 
    Note that for any $m<\alpha$, $f_0$ being an $\alpha$-H\"older function implies that $f_0$ is also an $m$-H\"older function. Therefore, the analysis in this work can be easily adapted to the case where the relation between $m$ and $\alpha$ is unknown by simply replacing $\alpha$ with $\min\{m,\alpha\}$. 
    When $\bs = \zero$, $\psi_{\bk\bs}$ closely resembles the term ${\varphi_h(\bx-\bx_i)}/{\sum_i\varphi_h(\bx - \bx_i)}$
    appearing in the summation of the Nadaraya-Watson estimator \eqref{eqn:NW_estimator}. When $\bs\neq \zero$, $\psi_{\bk\bs}$ is able to capture higher-order curvature information from the underlying regression function due to its polynomial structure.

   Finally, with the above ingredients, we define the kernel mixture of polynomials model by 
  $\calM = \left\{p_{f,\sigma}({\bx},y):f\in\bigcup_{K=1}^\infty\calF_K,\sigma\in[\lsigma,\usigma]\right\}$, a class of distributions indexed by the regression function $f$ and the standard deviation of the noise $\sigma$,
  where $p_{f,\sigma}({\bx},y)=\phi_\sigma\left({y-f({\bx})}\right)p_{\bx}({\bx})$, $f$ lies in the union of the function classes $\calF_K=\bigcup_{Kh\in[\llh,\uh]}\calF_K(h)$ for some constants $\llh,\uh$ with $1<\llh<\uh<\infty$, and
  \begin{align}
  \label{eqn:F_K}
  \calF_K(h)=\left\{\sum_{{\bk}\in[K]^p}\sum_{\bs:|\bs|=0}^m \xi_{\bk\bs}\psi_{\bk\bs}(\bx):
  {\bmu}_{\bk}\in\calX_K({\bk}),
  \max_{|\bs|=0,\ldots,m}|\xi_{\bk\bs}|\leq B,{\bk}\in[K]^p\right\}&.
  \end{align}
  The parameters $\sigma$, $K$, $(\bmu_\bk:\bk\in[K]^p)$, and $(\xi_{\bk\bs}:\bk\in[K]^p,|\bs|=0,1,\ldots,m)$ are to be assigned a hierarchical prior in Section \ref{sub:prior_specification}. In other words, the kernel mixture of polynomials model parametrizes the regression function $f$ through
  \begin{align}
   \label{eqn:kernel_mixture_polynomial}
  f({\bx})=\sum_{{\bk}\in[K]^p}\sum_{\bs:|\bs|=0}^m\xi_{\bk\bs}\psi_{\bk\bs}(\bx),
   \end{align}
  where $\{\psi_{\bk\bs}(\bx)\}_{\bk\bs}$ serves as certain basis functions that mimic the behavior of $f$ locally around $\bmu_\bk\in\calX_K(\bk)$, and $\xi_{\bk\bs}$ is the coefficient or the amplitude of the corresponding basis function $\psi_{\bk\bs}(\bx)$.
  \begin{remark}
  In order that the centers $\bmu_\bk$'s in the kernel mixture weights are helpful to detecting both the local and global behavior of $f$, they need to be relatively spread. 
  This also explains why in the first step above, we require that $\bmu_k$'s lie in the disjoint blocks $\in\calX_K(\bk)$'s. 
  In addition, such a partition restriction avoids the label-switching phenomenon in the Markov chain Monte Carlo sampler for posterior inference, which could potentially affect the mixing of the Markov chain \citep{celeux2000computational,jasra2005markov}. Alternatively, repulsive priors \citep{affandi2013approximate,xu2016bayesian,xie2017bayesian} can be incorporated to gain well-spread kernel centers, but this could still  cause the label-switching issue.
  \end{remark}

  \begin{remark}
  The kernel bandwidth parameter $h$, which is typically endowed with a prior in the literature of Bayesian kernel methods \citep{ghosal2007convergence,shen2013adaptive}, also plays a key role in establishing the convergence properties of the kernel mixture of polynomials. 
  In the current context, $h$ can be arbitrarily close to zero asymptotically, since we require that $K$ ranges over all positive integers and that $Kh$ stays bounded as $K\to\infty$. Therefore, $\bigcup_{K=1}^\infty\calF_K$ is rich enough to provide good approximation to arbitrary $f_0\in\mathfrak{C}^{\alpha,B}(\calX)$. 
  \end{remark}

  \subsection{Prior specification} 
  \label{sub:prior_specification}
   We define a prior distribution $\Pi$ for $(f,\sigma)$ through \eqref{eqn:kernel_mixture_polynomial} by imposing the hierarchical priors on the parameters $({\bmu}_{\bk}:{\bk}\in[K]^p)$, $(\xi_{\bk\bs}: \bk\in[K]^p, |\bs|=0,1,\ldots,m)$, $\sigma$, $h$, and $K$ as follows:
  \begin{itemize}
    \item The standard deviation $\sigma$ of the noises $(e_i)_{i=1}^n$ follows $\pi_\sigma$ that is continuous and non-vanishing on $[\lsigma{},\usigma{}]$, independent of the remaining parameters. 
    \item The prior for $K$ satisfies the following condition:
    \begin{align}\label{eqn:prior_rate_pi_K}
    \exp\left[-b_0x^p\left(\log x^p\right)^{r_0}\right]\lesssim\Pi(K\geq x)\lesssim \exp\left[-b_1x^p\left(\log x^p\right)^{r_0}\right]
    \end{align}
    for some constants $b_0,b_1>0$, $b_1\leq b_0$, and $r_0\geq0$. 
  \end{itemize}
  Given $K$, the conditional prior for the rest of the parameters are given as follows:
  \begin{itemize}
    \item The kernel centers $\bmu_\bk$'s are sampled as ${\bmu}_{\bk}={\bmu}_{\bk}^\star+\widetilde{{\bmu}}_{\bk}/(2K)$, where $\widetilde{{\bmu}}_{\bk}$ independently follows $\pi_{\bmu}$ for each ${\bk}\in[K]^p$ for some continuous non-vanishing density $\pi_\bmu$ on the hypercube $[-1,1]^p$.
    For example, $\pi_\bmu$ can be taken as the uniform distribution on $[-1,1]^p$. 
    \item The coefficient $\xi_{\bk\bs}$ for $\psi_{\bk\bs}(\bx)$ follow $\pi_\xi$ independently for each $\bs\in[m]^p$ with $|\bs|=0,1,\ldots,m$ given $K$, where $\pi_\xi$ is a continuous non-vanishing density on $[-B,B]$. For example, one can take $\pi_\xi$ to be the normal distribution truncated on $[-B,B]$, or the uniform distribution on $[-B,B]$.
    \item The bandwidth $h$ follows a non-vanishing density $\pi(h\mid K)$ supported on $[\llh/K,\uh/K]$ given $K$, where $1<\llh<\uh$. The uniform distribution on $[\llh/K,\uh/K]$ satisfies the condition.
  \end{itemize}
  \begin{remark}
  \cite{kruijer2010adaptive} adopts the same tail probability condition \eqref{eqn:prior_rate_pi_K} for the number of support points in the context of nonparametric density estimation. Special cases of \eqref{eqn:prior_rate_pi_K} include the geometric distribution when $r_0=0$, and the Poisson distribution when $r_0=1$. In Section \ref{sec:theoretical_properties_of_the_bayesian_kernel_regression}, we show that both $r_0=0$ and $r_0=1$ yield the same rate of contraction, but any $r_0>1$ ({i.e.}, thinner tail) leads to a slower rate of contraction.
  \end{remark}

  \section{Convergence properties of the kernel mixture of polynomials regression} 
  \label{sec:theoretical_properties_of_the_bayesian_kernel_regression}
  In this section, we establish the convergence results of the kernel mixture of polynomials for nonparametric regression, and obtain a frequentist sieve maximum likelihood estimator with the corresponding convergence rate. 
  For nonparametric regression problems,  
  when the true regression function $f_0$ is in $\mathfrak{C}^{\alpha,B}(\calX)$, $p_{\bx}({\bx})=1$, and $e_i\sim\mathrm{N}(0,1)$, $i=1,\ldots,n$, it has been shown that 
  the {minimax rate of convergence} for any estimator with respect to the $L_2$-norm
  is $n^{-{\alpha}/{(2\alpha+p)}}$ \citep{stone1982optimal,gyorfi2006distribution}. The optimal rate of contraction cannot be faster than the minimax rate of convergence. 
  {Theorem \ref{thm:contraction_rate}} below, which is one of the main results of this section, asserts that the rate of contraction with respect to the $L_2(\mathbb{P}_\bx)$-topology
  is minimax-optimal up to a logarithmic factor. Furthermore, the rate of contraction is adaptive to the smoothness level $\alpha$ of the underlying true $f_0$. 
    
  \begin{theorem}[Rate of contraction]\label{thm:contraction_rate}
   Assume that $f_0$ is in the $\alpha$-H\"older function class $\mathfrak{C}^{\alpha,B}(\calX)$ with envelope $B$. Suppose $\Pi$ is the prior constructed in Section \ref{sub:prior_specification}. Then for some large constant $M>0$, it holds that 
   \begin{align}\label{eqn:contraction_rate_L2}
   \Pi(\left\|f-f_0\right\|_{L_2(\mathbb{P}_\bx)}^2>M\eps_n^2\mid\calD_n)\to 0
   \end{align}
   in $\mathbb{P}_0$-probability, where $\calD_n$ denotes the data $(\bx_i,y_i)_{i = 1}^n$, $\eps_n = n^{-{\alpha}/{(2\alpha+p)}}(\log n)^{t/2}$, and $t>{2\alpha}\max(r_0,1)/{(2\alpha+p)}+\max\left(0,{1-r_0}\right)$. 
   \end{theorem}


   We sketch the proof below and defer the details to the Supplementary Material. Write the posterior distribution $\Pi(\cdot\mid\calD_n)$ as follows:
   \[
   \Pi(\calA\mid \calD_n) = \frac{\int_\calA \exp[\ell_n(f,\sigma) - \ell_n(f_0,\sigma_0)]\Pi(\mathrm{d}f\mathrm{d}\sigma)}{\int \exp[\ell_n(f,\sigma) - \ell_n(f_0,\sigma_0)]\Pi(\mathrm{d}f\mathrm{d}\sigma)} := \frac{\mathfrak{N}_n(\calA)}{\mathfrak{D}_n}
   \]
   for any measurable set $\calA$, where $\ell_n(f,\sigma)$ is the log-likelihood function
   $\ell_n(f, \sigma) = \sum_{i = 1}^n \log p_{f,\sigma}(\bx_i,y_i)$.
   To prove the rates of contraction results in Theorem \ref{thm:contraction_rate}, it suffices to verify a set of sufficient conditions, one variation of the prior-concentration-and-testing framework originally presented by \cite{ghosal2000convergence}. 
   This framework, along with other variations, have been widely applied not only to nonparametric regression, but also to density estimation \citep{ghosal2001entropies,ghosal2007posterior,kruijer2010adaptive,shen2013adaptive} and high-dimensional statistics \citep{castillo2012needles,bhattacharya2015dirichlet,rovckova2018bayesian,pati2014posterior,gao2015rate}.

   For the purpose of illustrating the proof framework, we introduce the notions of covering number and metric entropy. 
   For a metric space $(\calF, d)$, for any $\epsilon>0$, the $\epsilon$-covering number of $(\calF,d)$, denoted by $\calN(\epsilon,\calF,d)$, is defined to be the minimum number of $\epsilon$-balls of the form $\{g\in\calF:d(f,g)<\epsilon\}$ that are needed to cover $\calF$. The $\epsilon$-bracketing number of $(\calF,d)$, denoted by $\calN_{[\bdot]}(\epsilon,\calF,d)$, is defined to be the minimum number of brackets of the form $[l_i,u_i]$ that are needed to cover $\calF$ such that $l_i,u_i\in\calF$ and $d(l_i,u_i)<\epsilon$. We refer to $\log\calN(\epsilon,\calF,d)$ as the metric entropy, and $\log\calN_{[\bdot]}(\epsilon,\calF,d)$ as the bracketing (metric) entropy. The bracketing integral $\int_0^\epsilon\sqrt{\log\calN_{[\bdot]}(u,\calF,d)}\mathrm{d}u$ is denoted by $J_{[\bdot]}(\epsilon,\calF,d)$.

   We now specify the prior-concentration-and-testing framework, which are two sufficient conditions for \eqref{eqn:contraction_rate_L2}. 
   The first one is a prior concentration condition. Namely, there exists another sequence $(\leps_n)_{n=1}^\infty$ with $\leps_n\leq \eps_n$, such that
     \[
     \Pi(p\in B_{\mathrm{KL}}(p_0,\leps_n))\geq\exp(-n\leps_n^2),
     \]
     where $B_{\mathrm{KL}}(p_0,\eps)$ is the Kullback-Leibler ball of radius $\eps$ centered at $p_0$ defined as follows:
     \[
     B_{\mathrm{KL}}(p_0,\epsilon)=\left\{p_{f,\sigma}:D_{\mathrm{KL}}(p_0||p_{f,\sigma})<\epsilon^2,\mathbb{E}_0\left[\log \frac{p_0(\bx,y)}{p_{f,\sigma}(\bx,y)}\right]^2<\epsilon^2\right\}.
     \]
     Secondly, we  require a summability condition: For the sequence $(\leps_n)_{n = 1}^{\infty}$ above, there exists
     a sequence of sub-classes of densities $(\calM_n)_{n=1}^\infty$, $\calM_n\subset\calM$ (recall the definition of $\calM$ in Section \ref{sub:bayesian_kernel_smoother}), and for each $\calM_n$ a partition $(\calM_{nm})_{m=1}^\infty$ with $\calM_n=\bigcup_{m=1}^\infty\calM_{nm}$, such that 
     $\Pi(p\in \calM_n^c)\leq \exp(-4n\leps^2_n)$ and 
     \[
     \exp(-n\eps_n^2)\sum_{m=1}^\infty\sqrt{\calN(\eps_n,\calM_{nm},H)}\sqrt{\Pi(p\in\calM_{nm})}\to 0.
     \]
     The sequence of sub-classes of densities $(\calM_n)_{n=1}^\infty$ is referred to as sieves in the literature \citep{shen1994convergence}.

   The prior concentration condition plays a fundamental role in proving convergence of posterior distributions, as it guarantees that with probability tending to one, the denominator $\mathfrak{D}_n$ appearing in the posterior distribution $\Pi(\cdot\mid\calD_n)$ does not decay super-exponentially (see, for example, Lemma 8.10 in \citealp{ghosal2017fundamentals}):
   \[
   \mathbb{P}_0\left(\mathfrak{D}_n\geq\Pi(p \in B_{\mathrm{KL}}(p_0, \leps_n))\exp[-(1 + c)n\leps_n^2]\right)\to 1
   \]
   for any positive $c > 0$. 
   Under the current setup, the prior concentration condition largely depends on how well the 
   kernel mixture of polynomials is able to approximate any $\alpha$-H\"older function. To this end, we introduce the following lemma to demonstrate the approximation power of the kernel mixture of polynomials.
   \begin{lemma}[Approximation lemma]\label{lemma:prior_thickness_lemma}
  Assume that $f_0$ is in the $\alpha$-H\"older function class $\mathfrak{C}^{\alpha,B}(\calX)$ with envelope $B$. 
  Let $f$ be of the form \eqref{eqn:kernel_mixture_polynomial} and ${\bmu}_{\bk}\in\calX_K({\bk})$. Then there exists some constant $C_1$ such that for sufficiently small $\epsilon$ the following holds whenever $K\geq \epsilon^{-1/\alpha}$,
  \begin{align}
    B_K^\star:=&\left\{f:\max_{\bk\in[K]^p,|\bs|=0,1,\ldots,\lceil\alpha-1\rceil}\left|\xi_{\bk\bs}-\frac{D^\bs f_0({\bmu}_{\bk}^\star)}{s_1!\ldots s_p!}\right|\leq\epsilon, {\bk}\in[K]^p\right\}
    \nonumber\\
    \subset&\left\{f:\left\|f-f_0\right\|_{L_2(\mathbb{P}_\bx)}^2
    <C_1\epsilon^2\right\}.
    \nonumber
    \end{align}
   \end{lemma}
The intuition of Lemma \ref{lemma:prior_thickness_lemma} is that for each $\bmu_\bk^\star$, the function $f_0$ can be well approximated by the Taylor polynomials of degree $m$ locally around $\bmu_\bk^\star$. In fact, the coefficients $\xi_{\bk\bs}$ are selected to be sufficiently close to the Taylor polynomial coefficients. Then the local approximation effect around $\bmu_\bk^\star$ for each $\bk\in[K]^p$ is accumulated through the kernel mixture weights $(w_\bk(\bx):\bk\in[K]^p)$. Lemma \ref{lemma:prior_thickness_lemma} may be of independent interest for numerical function approximation as well. 

Moving forward to the summability condition, it states that there exists a subset $\calM_n$ of the entire model $\calM$ that occupies most of the prior probability, and at the same time can be covered by the union of a collection of blocks $(\calM_{nm})_{m\geq1}$ with low model complexity (metric entropies). Therefore, it is desired that useful metric entropy bounds can be obtained for the kernel mixture of polynomials model. The following proposition, which is one of the major technical contributions of this paper, directly tackles this issue.
  \begin{proposition}[Metric entropy bound]\label{lemma:metric_entropy_bound}
  There exists some constant $c_2>0$, such that for sufficiently small $\epsilon>0$ and any $r\in[1,\infty)$, 
  \begin{align}
  \log\calN_{[\bdot]}(2\epsilon, \calF_K, \|\cdot\|_{L_r(\mathbb{P}_\bx)})
  &\leq\log\calN(\epsilon, \calF_K, \|\cdot\|_\infty)
  \leq  c_2K^p\left(\log\frac{1}{\epsilon}\right)\nonumber.
  \end{align}
  \end{proposition}

   Besides the rate of contraction, which is a frequentist large sample evaluation of the full posterior distribution, we also obtain a frequentist sieve maximum likelihood estimator with a convergence rate as a result of the metric entropy bounds. This convergence rate is also minimax optimal up to a logarithmic factor. Interestingly, this rate is tighter than the rate of contraction of the full posterior, but the price we pay for the rate improvement is that the construction of the sieve depends on the smoothness level $\alpha$ and the rate is non-adaptive. 
  \begin{theorem}\label{thm:convergence_rate_sieveMLE}
  Assume that $f_0$ is in the $\alpha$-H\"older function class $\mathfrak{C}^{\alpha,B}(\calX)$ with envelope $B$. 
  Consider the sieve maximum likelihood estimator $\widehat{f}_K({x})$ defined by
  \begin{align}
  \widehat{f}_K({\bx})=\argmax_{f\in\calG_K}\sum_{i=1}^n\log \phi_{\sigma_0}(y_i-f({\bx}_i)),\nonumber
  \end{align}
  where
  \begin{align}
  \calG_{K}=\left\{\sum_{{\bk}\in[K]^p}
  \sum_{\bs:|\bs|=0}^m\xi_{\bk\bs}\psi_{\bk\bs}(\bx):Kh\in\left[\llh,\uh\right],
  {\bmu}_{\bk}\in\overline{\calX_K({\bk})},\max_{\bk\in[K]^p,|\bs|=0,\ldots,m}|\xi_{\bk\bs}|\leq B\right\}\nonumber.
  \end{align}
  If $K_n=\lceil\left(n\middle/\log n\right)^{{1}/{(2\alpha+p)}}\rceil$, then 
  \[\lim_{n\to\infty}\mathbb{P}_0\left(\left\|f_0-\widehat{f}_{K_n}\right\|_{L_2(\mathbb{P}_\bx)}\geq M({\log n}/{n})^{{\alpha}/{(2\alpha+p)}}\right)=0\]
  for some large constant $M>0$.
  \end{theorem}


\section{Application to the partial linear model} 
\label{sec:bayesian_kernel_smoothing_in_partial_linear_models}
The kernel mixture of polynomials model enjoys flexible extendability thanks to its nonparametric nature. 
In this section we present a semiparametric application: we use the kernel mixture of polynomials to model the nonparametric component in the partial linear model. The partial linear model is of the form $y_i=\bz_i\transpose{}\bbeta+\eta({\bx}_i)+e_i$, where $\bz_i,{\bx}_i$'s are design points, $\bbeta$ is the linear coefficient, $\eta:\calX\to\mathbb{R}$ is an unknown nonparametric function, and $e_i$'s are independent $\mathrm{N}(0,1)$ noises. 
In many applications, the estimation of the nonparametric component $\eta$ is of great interest. For example, in \cite{xuyanbo2016bayesian}, the parametric term $\bz_i\transpose{}\bbeta$ models the baseline disease progression and the nonparametric term $\eta({\bx}_i)$ models the individual-specific treatment effect deviations over time. When the regression coefficient $\bbeta$ is of more interest, the estimation of $\eta$ can still be critical since it could affect the inference of $\bbeta$. 
As will be seen later, we prove 
 the convergence results for both the nonparametric component $\eta$ (Theorem \ref{thm:nonparametric_rate}) and the parametric component $\bbeta$ (Theorem \ref{thm:asymptotic_normality}). 
Furthermore, as a consequence of the metric entropy result (Proposition \ref{lemma:metric_entropy_bound}), 
we obtain the Bernstein-von Mises limit of the marginal posterior distribution of  $\bbeta$. 

\subsection{Setup and prior specification} 
\label{sub:setup_and_prior_specification}
Let $\calX=[0,1]^p\subset\mathbb{R}^p$ be the design space of the nonparametric component, $\calZ\subset\mathbb{R}^q$ be the design space of the parametric component, and $p_{({\bx},\bz)}:\calX\times\calZ\to(0,\infty)$ be a continuous density function supported on $\calX\times\calZ$. 
We incorporate the partial linear model with the kernel mixture of polynomials prior for the nonparametric component $\eta$ through
$\calP=\left\{p_{\bbeta,\eta}({\bx},\bz,y):\bbeta\in\mathbb{R}^q,\eta\in\bigcup_{K=1}^\infty\calF_K\right\}$, a class of distributions indexed by the linear coefficient $\bbeta$ and the nonparametric component $\eta$,
where 
$p_{\bbeta,\eta}({\bx},\bz,y)=\phi(y-\bz\transpose{}\bbeta-\eta({\bx}))p_{{\bx},\bz}({\bx},\bz)$ and $\calF_K=\bigcup_{Kh\in\left[\llh,\uh\right]}\calF_K(h)$ with $\calF_K(h)$ given by \eqref{eqn:F_K}. 
We assume that the data $\calD_n=({\bx}_i,\bz_i,y_i)_{i=1}^n$ are independently sampled from $p_0({\bx},\bz,y)=\phi(y-\bz\transpose{}\bbeta_0-\eta_0({\bx}))p_{{\bx},\bz}({\bx},\bz)$ for some $\bbeta_0\in\mathbb{R}^q$ and some function $\eta_0\in\mathfrak{C}^{\alpha,L}(\calX)$. Several additional assumptions regarding the parametric component $\bz\transpose\bbeta$ are needed: The design space $\calZ\subset\mathbb{R}^q$ for $\bz$ is compact with $\sup_{\bz\in\calZ}\|\bz\|_1\leq\overline{B}$ for some $\overline{B} > 0$; The sampling distribution for $\bz$ satisfies $\mathbb{E}\bz=\zero$ and $\mathbb{E}\bz\bz\transpose{}$ being non-singular; The density of the design points $({\bx}_i,\bz_i)_{i=1}^n$ factorizes as  $p_{({\bx},\bz)}({\bx},\bz)=p_{\bx}({\bx})p_\bz(\bz)$, {i.e.}, ${\bx}$ and $\bz$ are independent.

For the prior specification, we assume $\eta$ follows the kernel mixture of polynomials prior $\Pi_\eta$ constructed in Section \ref{sub:prior_specification} with $\sigma=1$. 
For the parametric component $\bbeta$, we impose a standard Gaussian prior $\Pi_\bbeta=\mathrm{N}(0,\eye_q)$, independent of $\Pi_\eta$. The joint prior is denoted by $\Pi=\Pi_\eta\times\Pi_\bbeta$.

\subsection{Convergence results} 
\label{sub:posterior_contraction_rate_for_the_nonparametric_component}
We first tackle the convergence of the nonparametric component $\eta$. The following theorem not only addresses the rate of contraction of the marginal posterior of $\eta$, but also serves as one of the building blocks for proving the Bernstein-von Mises limit of the marginal posterior of $\bbeta$. The proofs of Theorem \ref{thm:nonparametric_rate} and Theorem \ref{thm:asymptotic_normality} are deferred to the Supplementary Material. 
  \begin{theorem}[Nonparametric rate]\label{thm:nonparametric_rate}
  Assume that $\eta_0$ is in the $\alpha$-H\"older function class $\mathfrak{C}^{\alpha,B}(\calX)$ with envelope $B$. 
  Under the setup and prior specification in Section \ref{sub:setup_and_prior_specification},
   \[\Pi(\|\eta-\eta_0\|_{L_2(\mathbb{P}_\bx)}^2>M\eps_n^2\mid\calD_n)\to 0\]
  in $\mathbb{P}_0$-probability for some large constant $M>0$, where $\eps_n = n^{-{\alpha}/{(2\alpha+p)}}(\log n)^{t/2}$, and $t>{2\alpha}\max(r_0,1)/{(2\alpha+p)}+\max\left(0,{1-r_0}\right)$.  
  \end{theorem} 

Now we turn to the convergence results for the parametric component. The focus is the asymptotic normality of the marginal posterior distribution of $\bbeta$, {i.e.}, the Bernstein-von Mises limit \citep{doob1949application}. 
To achieve this, we need the notion of the \emph{least favorable submodel} for semiparametric models \citep{bickel1998efficient}. For each fixed $\bbeta\in\mathbb{R}^q$, the {least favorable curve} $\eta^*_\bbeta$ is defined by the minimizer of the Kullback-Leibler divergence over all $\eta$:
    $\eta^*_\bbeta(x)=\arginf_{\eta\in\calF}D_{\mathrm{KL}}(p_0||p_{\bbeta,\eta})$. 
Under the assumptions $\mathbb{E}\bz = \zero$ and $p_{(\bx,\bz)}(\bx,\bz) = p_\bx(\bx)p_\bz(\bz)$, for each $\bbeta$, it can be shown that
  $\eta^*_\bbeta({\bx})$ coincides with $\eta_0(\bx)$.
  The {least favorable submodel} is defined to be $\{p_{\bbeta,\eta^*_\bbeta}:\bbeta\in\mathbb{R}^q\}$, which in turn coincides with $\{p_{\bbeta,\eta_0}:\bbeta\in\mathbb{R}^q\}$ in our context. 

  \begin{theorem}\label{thm:asymptotic_normality}
  Assume that $\eta_0$ is in the $\alpha$-H\"older function class $\mathfrak{C}^{\alpha,B}(\calX)$ with envelope $B$. Under the setup and prior specification in Section \ref{sub:setup_and_prior_specification}, if $\alpha>p/2$, then
  \begin{align}
  \sup_{F}\left|\Pi\left(\sqrt{n}(\bbeta-\bbeta_0)\in F\mid\calD_n\right)-\Phi(F\mid\bDelta_n,(E\bz\bz\transpose{})^{-1})\right|\to 0\nonumber
  \end{align}
  in $\mathbb{P}_0$-probability, where $\Phi(\bdot\mid\bDelta_n,(\mathbb{E}\bz\bz\transpose{})^{-1})$ is the $\mathrm{N}(\bDelta_n,(\mathbb{E}\bz\bz\transpose{})^{-1})$ probability measure and
  \[\bDelta_n=\frac{1}{n}\sum_{i=1}^n(\mathbb{E}\bz\bz\transpose{})^{-1}\bz\left[y_i-\eta_0({\bx}_i)-\bz_i\transpose{}\bbeta_0\right].\]
  \end{theorem}
  The proof of Theorem \ref{thm:asymptotic_normality} is based on verifying a set of sufficient conditions in \cite{yang2015semiparametric}, which are provided in the Supplementary Material. However, we remark that the metric entropy results (Proposition \ref{lemma:metric_entropy_bound}) in Section \ref{sec:theoretical_properties_of_the_bayesian_kernel_regression} and the  rate of contraction for $\eta$ (Theorem \ref{thm:nonparametric_rate}) are also of fundamental interest in the verification process. 

  \section{Numerical studies} 
\label{sec:numerical_studies}
We perform numerical studies to evaluate the kernel mixture of polynomials for nonparametric regression and the partial linear model. 
The posterior inference for all examples is carried out by a Markov chain Monte Carlo (MCMC) sampler with the number of burn-in  iterations being $1000$. Then we collect $1000$ post-burn-in MCMC samples for posterior analysis. 
To determine $K$, we  collect posterior samples for each fixed $K\in\{K_{\min},\ldots,K_{\max}\}$ and find the optimal $K$ by minimizing the {deviance information criterion} \citep{gelman2014bayesian} over $K$.  Numerical evidence shows that all Markov chains converge within $1000$ iterations. The kernel we use throughout is the bump kernel $\varphi(x)=\exp[-1/(1-x^2)]\mathbf{1}(|x|<1)$. 

\subsection{A synthetic example for nonparametric regression} 
\label{sub:a_synthetic_example_for_nonparametric_regression}
We first consider a synthetic example for nonparametric regression. Following \cite{knapik2011bayesian} and \cite{yoo2016supremum}, we consider the true function to be
$f_0(x)=\sqrt{2}\sum_{s=1}^\infty s^{-{3}/{2}}\sin(s)\cos[(s-{1}/{2})\pi x]$,
which has smoothness level $\alpha=1$. We generate $n=1000$ observations $(x_i,y_i)_{i=1}^n$ given by $y_i=f_0(x_i)+e_i$, where the design points $(x_i)_{i=1}^n$ are independently and uniformly sampled over $\calX=[0,1]$, and $(e_i)_{i=1}^n$ are independent $\mathrm{N}(0,0.2^2)$ noises. For the prior specification, we let $\pi_\mu=\mathrm{Unif}(-1,1)$, $\pi_\beta=\mathrm{N}(0,10^2)$, $\pi_\xi=\mathrm{N}(0,10^2)\mathbf{1}(|\xi|\leq 50)$, and $Kh\sim\mathrm{Unif}(\llh,\uh)$. The hyperparameters are set as $\llh = 1.2,\uh = 2$, $B=50$, $K = 15$, and $m = 2$. The range of $K$ is set to be $\{6,7,\ldots,15\}$. 

\begin{figure}[!ht] 
\centerline{\includegraphics[width=.9\textwidth]{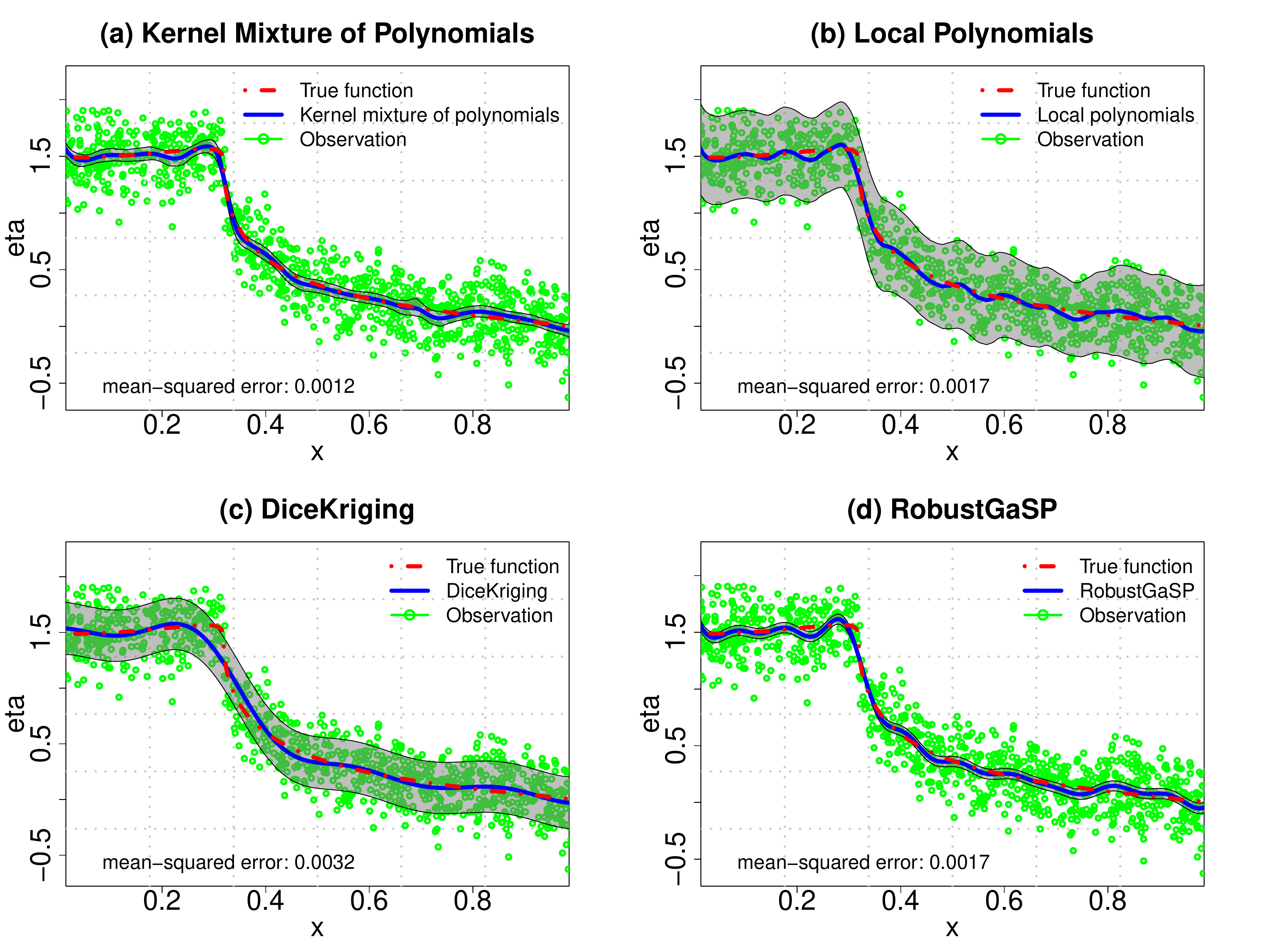}}
\caption{Synthetic example for nonparametric regression. Shaded regions are the point-wise $95\%$-credible/confidence intervals. The solid lines are point-wise posterior means/point estimators of $f$ and the dot-dashed lines are the ground true $f_0$. Scatter points are the observations. }
\label{Fig:simulation_volterra_function}
\end{figure}

For comparison we consider $3$ competitors for estimating $f_0$: the local polynomial regression \citep{fan1996local}, implemented in the \verb|locpol| package \citep{cabrera2012locpol}, DiceKriging method \citep{roustant2012dicekriging}, and the robust Gaussian stochastic process emulation (RobustGaSP, \citealp{Gu2017Robust}), implemented in the \verb|RobustGaSP| package \citep{Gu2016RobustGaSP}. 
The point-wise posterior means and $95\%$-credible/confidence intervals for $f(x)$ using the $4$ nonparametric regression methods are plotted in Figure \ref{Fig:simulation_volterra_function}, respectively. We also compute the mean-squared errors of the posterior means of the four methods, where the ground true $f_0$ is evaluated at $1000$ equidistant design points. In terms of accuracy measured by the mean-squared errors (marked in the bottom-left corner of each panel), the kernel mixture of polynomials performs better than DiceKriging and similarly to the local polynomials and RobustGaSP.

\begin{figure}[!ht] 
\centerline{\includegraphics[width=1.02\textwidth]{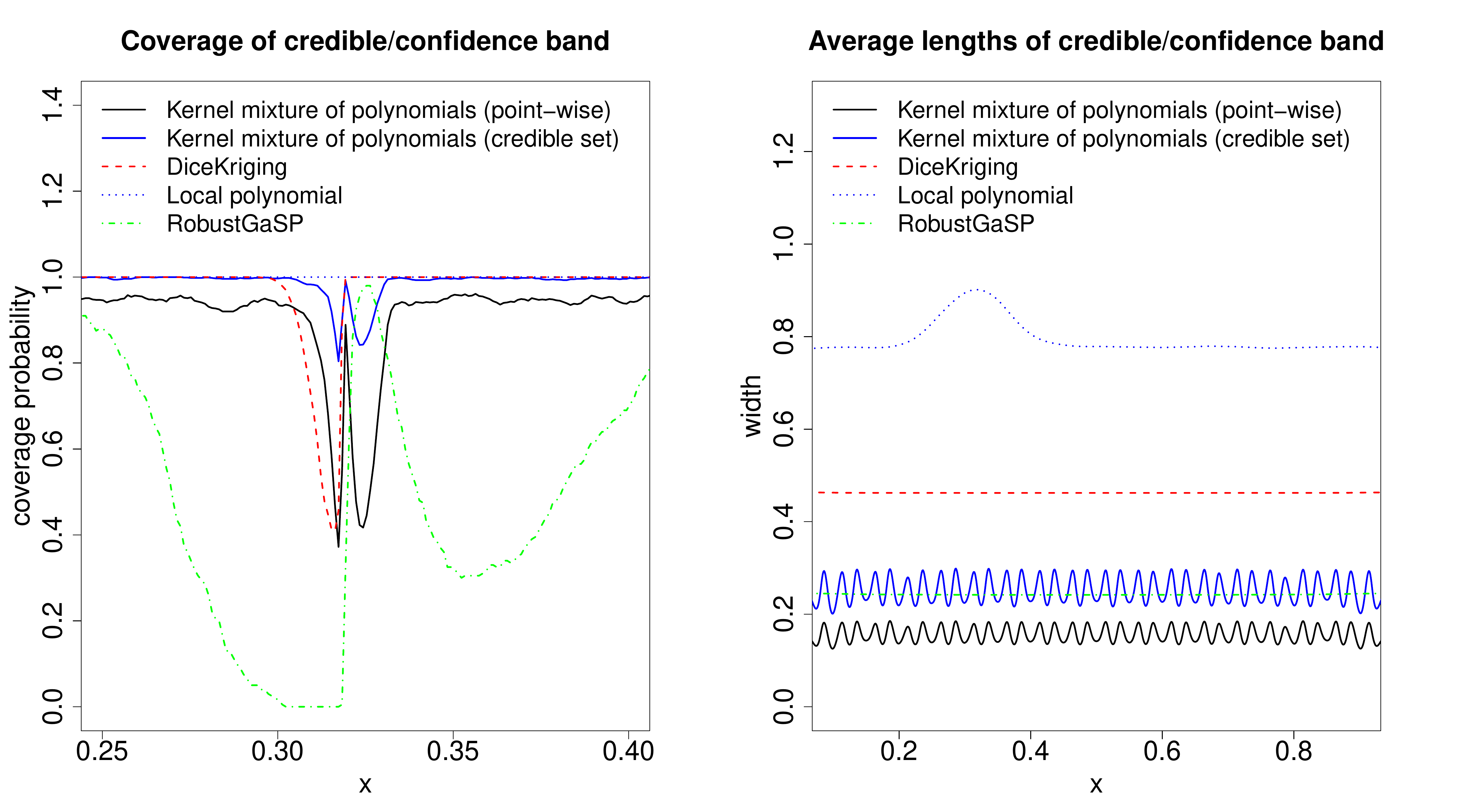}}
\caption{Coverage of the credible band/confidence band/credible set and the average lengths of the credible band/confidence band/credible set for nonparametric regression under $1000$ repeated experiments.}
\label{Fig:coverage_width_CI_regression}
\end{figure}

We further investigate the uncertainty quantification of the kernel mixture of polynomials and the $3$ competitors by performing $1000$ repeated experiments. Namely, we simulate $1000$ synthetic datasets using the same setup as above, perform the corresponding inference for each dataset, and then 
compute the average coverage and average lengths of the point-wise credible/confidence band.  We are particularly interested in the behavior of coverage of the point-wise credible/confidence band for $x\in[0.30, 0.35]$, where the true regression function $f_0$ exhibits a bump. 
As shown in the left panel of Figure \ref{Fig:coverage_width_CI_regression}, DiceKriging and the local polynomial regression provide the best coverage, but the widths of the point-wise confidence band are much larger (see the right panel of Figure \ref{Fig:coverage_width_CI_regression}). RobustGaSP provides the worst coverage among the $4$ approaches, and the corresponding confidence band is still wider than that of the kernel mixture of polynomials. In contrast, the proposed kernel mixture of polynomial regression provides better credible band coverage for $x\in[0.30, 0.35]$ than RobustGaSP. 
In addition, the average length of the point-wise credible band of the kernel mixture of polynomials is narrower than those of the other $3$ competitors. However, we also notice that this feature is at the cost of worse coverage than the local polynomial regression and DiceKriging. Namely, the narrower point-wise confidence band is over-confident near the bump of $f_0$ than DiceKriging and the local polynomial regression. 
One potential reason for this phenomenon could be that the model over-smooths the true regression function $f_0$ near the bump based on a random sample of limited size, but $f_0$ is much smoother elsewhere \citep{yoo2016supremum}. 

As pointed out by one of the referees, for the kernel mixture of polynomials, the bad coverage behavior of the point-wise credible band near the bump of $f_0$ may be alleviated by using the credible set for the entire function. Following  \cite{szabó2015}, we consider the following $L_2$-credible set for $f$: 
\[
\hat C_n(\gamma_n) = \{f: \|f - f_0\|_{L_2(\mathbb{P}_\bx)} \leq \gamma_n\},
\]
where $\gamma_n>0$ is the credible radius such that $\Pi(f \in \hat C_n(\gamma_n)\mid\calD_n) = 95\%$. Numerically, the credible set $C_n(\gamma_n)$ can be computed as follows: let $(f_t^{\mathrm{mc}})_{t = 1}^{n_{\mathrm{mc}}}$ be $n_{\mathrm{mc}}$ posterior samples of $f$  drawn from the MCMC. Let $\gamma_n$ be the $95\%$-quantile of $(\|f_t^{\mathrm{mc}} - \hat f\|_{L_2(\mathbb{P}_\bx)})_{t = 1}^{n_{\mathrm{mc}}}$, where $\hat{f}$ is the point-wise posterior mean function. Then the $L_2$-credible set $\hat C_n(\gamma_n)$ can be constructed by
\[
\left\{f: \min_{t:\|f_t^{\mathrm{mc}} - \hat f\|_{L_2(\mathbb{P}_\bx)}\leq \gamma_n}f_t^{\mathrm{mc}}(\bx)\leq f(\bx)\leq \max_{t:\|f_t^{\mathrm{mc}} - \hat f\|_{L_2(\mathbb{P}_\bx)}\leq \gamma_n}f_t^{\mathrm{mc}}(\bx)\text{ for all }\bx\in\calX\right\}.
\]
Using the aforementioned $1000$ replicates of synthetic datasets, we examine the corresponding point-wise coverage behavior of the $L_2$-credible sets of the entire function locally around the bump ($x\in[0.25, 0.4]$), compute the average lengths of these credible sets, and plot them in Figure \ref{Fig:coverage_width_CI_regression}. The $L_2$-credible sets exhibit better point-wise coverage behavior than the point-wise credible bands at the cost of slightly wider lengths. 
This is also in accordance with the observations in \cite{szabó2015}: the overall coverage is satisfactory, but the coverage near the bump  is slightly lower.
We also observe that the lengths of these $L_2$-credible sets are close to those obtained from RobustGaSP, but the coverage is better. 

\subsection{Computation comparison with Gaussian processes} 
\label{sec:comparison_with_rescaled_gaussian_processes}

One motivation of developing the kernel mixture of polynomials for Bayesian nonparametric regression is to tackle the computational bottleneck of  classical Gaussian processes. In this section we 
compare the computational costs of the kernel mixture of polynomials against Gaussian processes with $3$ different covariance functions: the squared-exponential covariance function: $$K(x,x';\psi) = \exp\left[-({x-x'})^2/{\psi}^2\right],$$
the Mat\'ern covariance function with roughness parameter $3/2$: $$K(x,x';\psi) = \left(1+{\sqrt{3}|x-x'|}/{\psi}\right)\exp\left(-{\sqrt{3}|x-x'|}/{\psi}\right),$$ and the Mat\'ern covariance function with roughness parameter $5/2$:
  $$K(x,x';\psi) = \left(1+{\sqrt{5}|x-x'|}/{\psi}+{5|x-x'|^2}/{(3\psi^2)}\right)\exp\left(-{\sqrt{5}|x-x'|}/{\psi}\right).$$
The parameter $\psi$ appearing in these covariance functions is referred to as the range parameter. Here we follow the suggestion by \cite{van2009adaptive} and impose an inverse-Gamma hyperprior on $\psi$, \emph{i.e.}, $\pi(\psi)\propto\psi^{-a_\psi-1}\exp(-b_\psi/\psi)$ for some $a_\psi,b_\psi\geq 2$. Such a hierarchical formulation for Gaussian processes with a hyperprior on the range parameter is also referred to as rescaled Gaussian processes in the literature \citep{van2007bayesian}. 

\begin{figure}[!ht] 
\centerline{\includegraphics[width=1\textwidth]{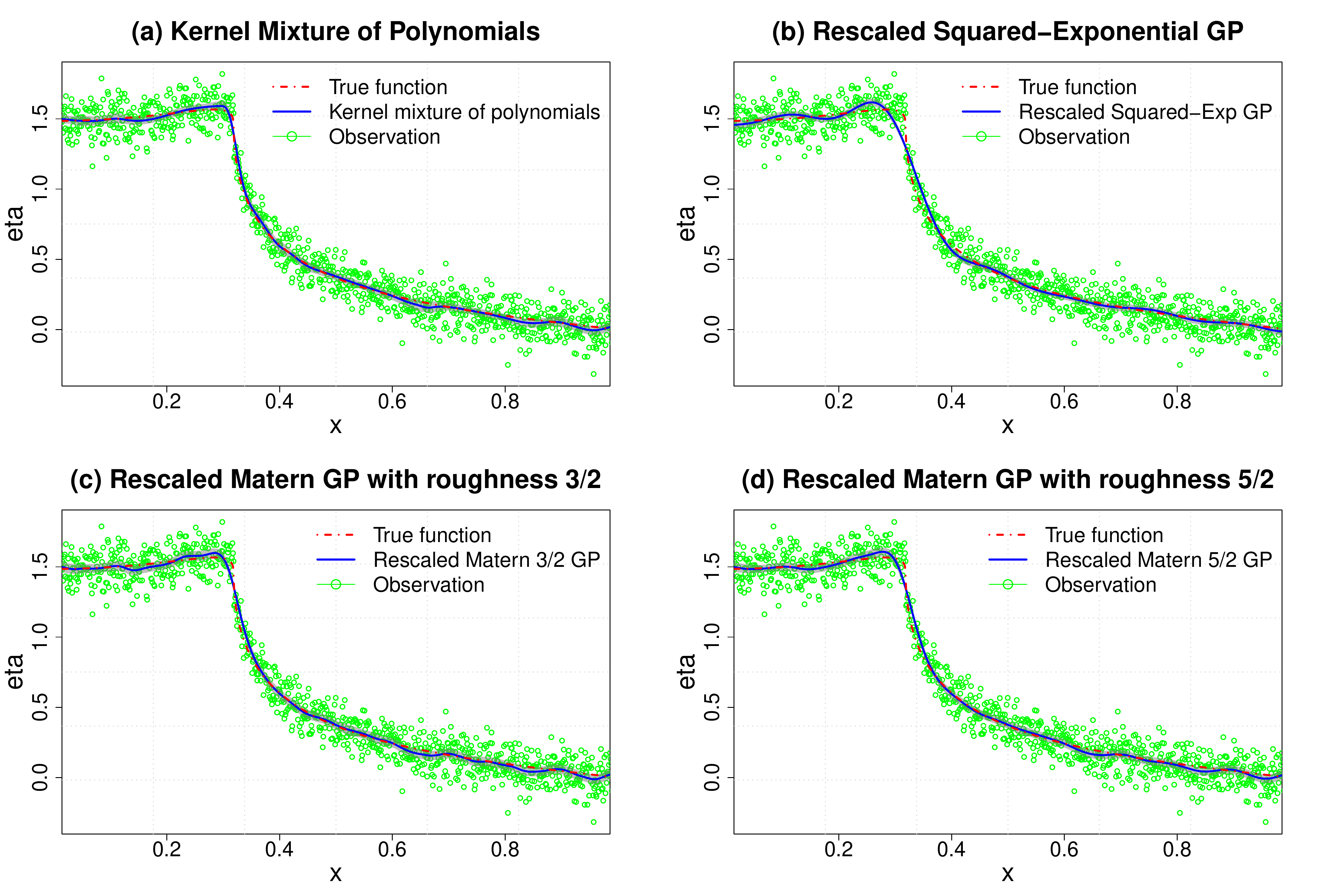}}
\caption{Comparison with rescaled Gaussian processes: Shaded regions are point-wise $95\%$-credible/confidence intervals. The solid lines are point-wise posterior means/point estimators of $f$ and the dot-dashed lines are the ground true $f_0$. Scatter points are the observations. }
\label{Fig:KEMPO_vs_RescaledGP}
\end{figure}

We adopt the same setup as in Subsection \ref{sub:a_synthetic_example_for_nonparametric_regression} except that $e_i\sim\mathrm{N}(0,0.1^2)$, and for each rescaled Gaussian process described above, posterior inference is carried out using classical Markov chain Monte Carlo, with $1000$ posterior samples collected after burn-in. The hyperparameter for $\pi_\psi$ is set to be $a_\psi = b_\psi = 2$. The posterior means and point-wise $95\%$ credible intervals for the kernel mixture of polynomials and 3 rescaled Gaussian processes are visualized in Figure \ref{Fig:KEMPO_vs_RescaledGP} (a), (b), (c), and (d), respectively. Furthermore, we examine the mean-squared errors of the posterior means as well as the computational costs of the 4 approaches. As shown in Table \ref{Tab:Comparison_rescaled_GP},  the proposed kernel mixture of polynomials not only yields a point estimate with smaller mean-squared error, but is also much more efficient than the rescaled Gaussian processes in terms of runtime expenses. The main reason is that when updating the range parameter $\psi$ in a single iteration of MCMC using Gaussian processes, the inverse of the covariance matrix needs to be re-computed, whereas the kernel mixture of polynomials does not suffer from such a computational bottleneck.

\begin{table}[!ht]
  \centering
  \caption{Comparison with  Gaussian processes: Mean-squared errors and computational costs}
    \begin{tabular}{ccccc}
    \hline\hline
    Method & mean-squared error & computational cost \\
    \hline
    Kernel mixture of polynomials & $3.447\times 10 ^{-4}$ & $562$s\\
    Squared-exponential Gaussian process & $1.266\times 10^{-3}$ & $198806$s\\
    Mat\'ern $3/2$ Gaussian process & $5.736\times 10^{-4}$ & $24219$s\\
    Mat\'ern $5/2$ Gaussian process & $6.918\times 10^{-4}$ & $20300$s\\
    \hline
    \hline
    \end{tabular}%
  \label{Tab:Comparison_rescaled_GP}
\end{table}%


\subsection{A synthetic example for the partial linear model} 
\label{sec:a_synthetic_example_for_partial_linear_model}
We consider a synthetic example to evaluate the performance of the kernel mixture of polynomials for the partial linear model.  We simulate $n = 500$ observations according to the model $y_i=\bz_i\transpose{}\bbeta_0+\eta_0(x_i)+e_i$, where $\bbeta_0$ is provided in Table \ref{Tab:simulation_beta}, $(e_i)_{i=1}^n$ are independent $\mathrm{N}(0,1)$ noises that are independent of $(x_i,\bz_i)_{i=1}^n$, and $\eta_0(x)=2.5\exp(-x)\sin(10\pi x)$. The nonparametric function $\eta_0$ is highly nonlinear and hence brings natural challenge to estimation. The design points $(\bz_i)_{i=1}^n$ for the linear component follow $\mathrm{Unif}([-1,1]^8)$ independently, and the design points $(x_i)_{i=1}^n$ for $\eta$ are independently sampled from $\mathrm{Unif}(0,1)$. The hyperparameters for the the kernel mixture of polynomials prior are set as follows: $\llh = 1.2,\uh = 2$, $B=100$, and $m = 3$. For the prior specification, we assume $\pi_\mu=\mathrm{Unif}(-1,1)$, $\pi_\bbeta=\mathrm{N}(\zero,10^2\eye_8)$, $\pi_\xi=\mathrm{N}(0,10^2)\mathbf{1}(|\xi|\leq 100)$, and $Kh\sim\mathrm{Unif}(\llh,\uh)$. We set the range of $K$ to be $\{6,7,\ldots,15\}$. 

\begin{figure}[!ht] 
\centerline{\includegraphics[width=.9\textwidth]{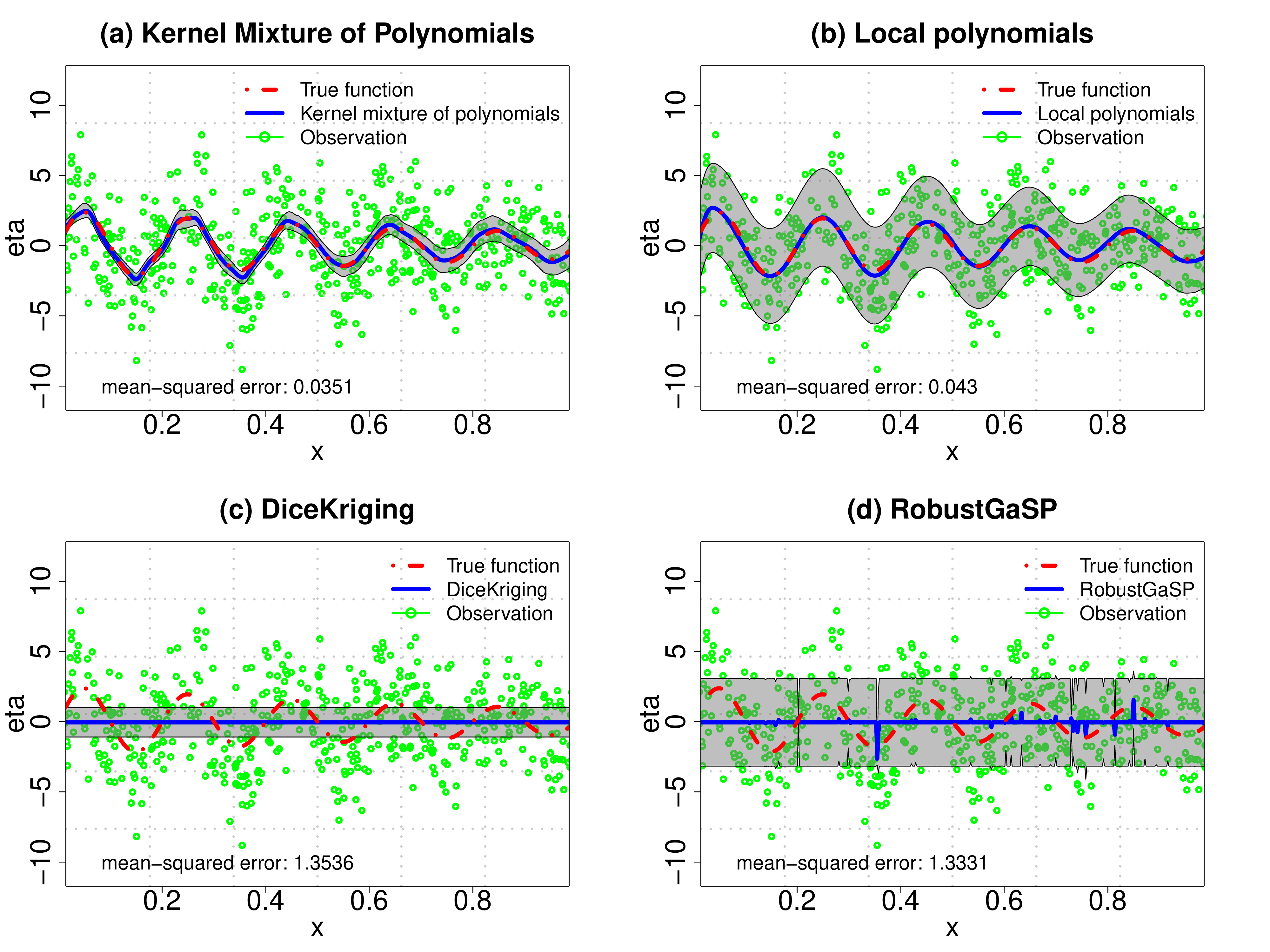}}
\caption{Synthetic example for the partial linear model. Shaded regions are point-wise $95\%$ credible/confidence intervals. The solid lines are point-wise posterior means/point estimators of $\eta$ and the dot-dashed lines are $\eta_0$. Scatter points are the observations.
Panels (b), (c), and (d) are computed using the discrepancy data $(x_i,y_i-\bz_i\transpose\widehat{\bbeta}_{\mathrm{LS}})_{i=1}^n$, where $\widehat{\bbeta}_{\mathrm{LS}}$ is the least-squared estimate of $\bbeta$. }
\label{Fig:simulation_nonparametric_component}
\end{figure}

For the parametric component, we compute 
the posterior means and the posterior $95\%$ credible intervals for $\bbeta$. For comparison, we calculate the least-squared estimate of $\bbeta$: 
$\widehat{\bbeta}_{\mathrm{LS}}=\argmin_{\bbeta\in\mathbb{R}^q}\sum_{i=1}^n(y_i-\bz_i\transpose{}\bbeta)^2.$
The comparison is provided in Table \ref{Tab:simulation_beta}, where the column ``KMP'' stands for the posterior means of $\bbeta$ under the kernel mixture of polynomials prior. From the posterior summary of $\bbeta$ we see that the underlying true $\bbeta_0$ lie in the posterior $95\%$-credible intervals of the kernel mixture of polynomials prior, and the corresponding posterior means outperform the least-squared estimate in terms of accuracy. 
\begin{table}[htbp]
  \centering
  \caption{Simulation example: Inference of $\bbeta$ with $\eta_0(x) = 2.4\exp(-x)\sin(10\pi x)$.}
    \begin{tabular}{ccccc}
    \hline\hline
    $\bbeta$ & $\bbeta_0$ & KMP  & $95\%$-credible intervals & $\widehat{\bbeta}_{LS}$ \\
    \hline
    $\beta_1$ & 1.0338 & 1.0579 & (0.8949, 1.2150) & 1.1976 \\
    $\beta_2$ & 0.1346 & 0.1003 & (-0.0560, 0.2516) & 0.1733 \\
    $\beta_3$ & 0.2854 & 0.3481 & (0.1832, 0.5090) & 0.4427 \\
    $\beta_4$ & 0.6675 & 0.6449 & (0.5007, 0.7887) & 0.8386 \\
    $\beta_5$ & 0.6732 & 0.7212 & (0.5630, 0.8838) & 0.7427 \\
    $\beta_6$ & 0.5293 & 0.5433 & (0.3971, 0.6866) & 0.6274 \\
    $\beta_7$ & -0.5073 & -0.4759 & (-0.6337, -0.3158) & -0.3464 \\
    $\beta_8$ & -3.3942 & -3.3031 & (-3.4492, -3.1450) & -3.5253 \\
    \hline
    \hline
    \end{tabular}%
  \label{Tab:simulation_beta}
\end{table}%
For the nonparametric component, we plot the point-wise posterior means of $\eta$ along with the point-wise $95\%$ credible intervals in Figure \ref{Fig:simulation_nonparametric_component} (a). The mean-squared error of the posterior means at $500$ equidistant points on $(0,1)$ is $0.0351$. Based on the least-squared estimate of $\bbeta$, we also consider the local polynomial regression, DiceKriging method, and RobustGaSP as three alternatives to estimate the nonparametric component.
The numerical comparisons are illustrated in Panels (b), (c), and (d) in Figure \ref{Fig:simulation_nonparametric_component}, respectively. The kernel mixture of polynomials outperforms the local polynomial in terms of both the mean-squared error. Both DiceKriging and RobustGaSP fail to detect the nonlinearity of $\eta_0$, giving rise to significantly larger mean-squared error of the point-wise posterior means. 

\begin{figure}[!ht] 
\centerline{\includegraphics[width=.9\textwidth]{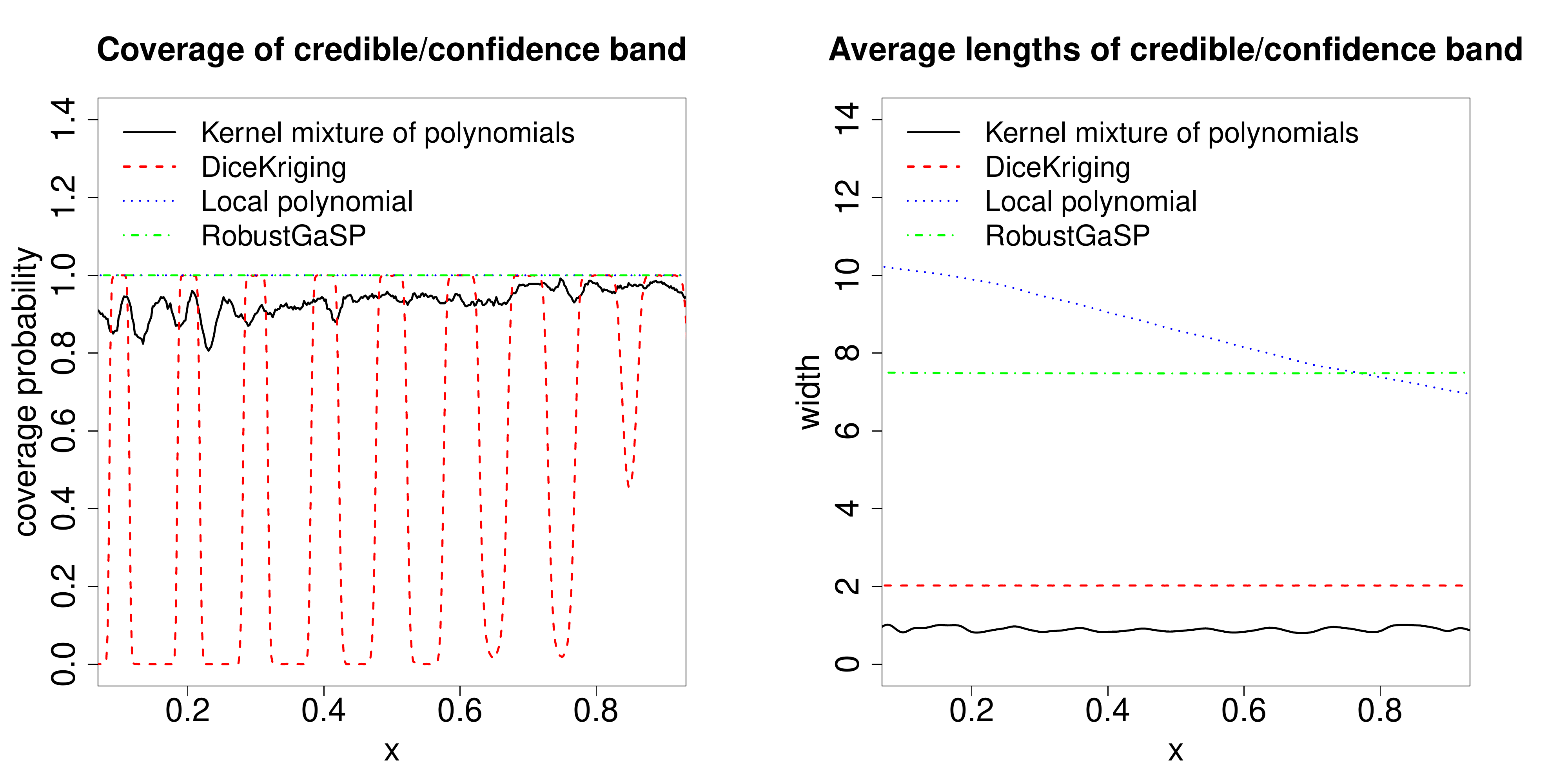}}
\caption{Coverage of point-wise credible/confidence band and average lengths of point-wise credible/confidence band for partial linear model under $1000$ repeated experiments.}
\label{Fig:coverage_width_CI_PLM}
\end{figure}

\begin{figure}[!ht] 
\centerline{\includegraphics[width=.9\textwidth]{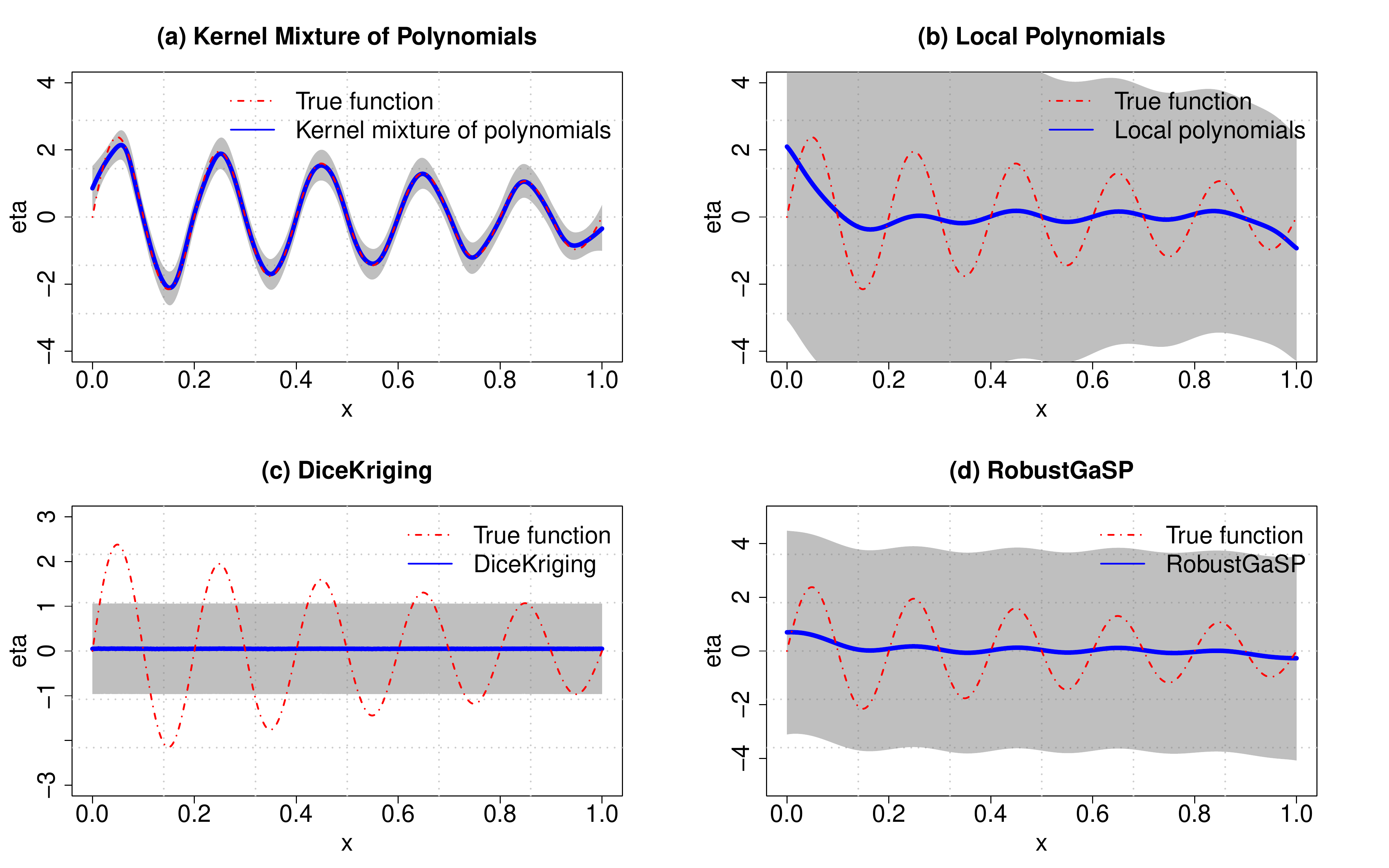}}
\caption{Averages of regression function estimates (blue solid lines) and the corresponding point-wise credible/confidence band (shaded regions) for partial linear model using $1000$ repeated  experiments.}
\label{Fig:PLM_repeat_experiments}
\end{figure}

To examine the uncertainty quantification of the $4$ approaches, we perform $1000$ repeated experiments. The coverage and average lengths of the point-wise credible/ confidence band are presented in Figure \ref{Fig:coverage_width_CI_PLM}. The averages of the regression function estimates and their point-wise credible/confidence bands based on $1000$ repeated experiments are computed and visualized in Figure \ref{Fig:PLM_repeat_experiments}. It can be seen that the kernel mixture of polynomials provides decent coverage to $\eta_0$ and the corresponding average lengths of the point-wise credible band are smaller than those using the other three approaches. In contrast, all the three competitors fail to detect the signal of $\eta_0$ after averaging from repeated  experiments, producing inappropriate point-wise confidence band (either too narrow or too wide). In particular, DiceKriging provides poor coverage of the point-wise confidence band, and the average widths of the point-wise confidence bands using the local polynomials and RobustGaSP are much wider than the magnitude of $\eta_0$.

\subsection{The partial linear model for the wage data} 
\label{sub:partial_linear_regression_for_the_wage_data}
We further analyze the cross-sectional data on wages \citep{wooldridge2015introductory}, a benchmark dataset for the partial linear model. This dataset is also available in the \verb|np| package \citep{hayfield2008nonparametric}. It consists of 526 observations with 24 variables and are taken from U.S. Current Population Survey for the year 1976. In particular, we are interested in modeling the hourly wage on the logarithm scale as the response with respect to 5 variables: years of education (educ), years of potential experience (exper), years with current employer, gender, and marital status. \cite{choi2015partially} and \cite{hayfield2008nonparametric} suggest the following form of the model:
\[
y_i = \beta_1z_i^{\mathrm{female}}+\beta_2z_i^{\mathrm{married}}+\beta_3z_i^{\mathrm{educ}}+\beta_4z_i^{\mathrm{tenure}}+\eta(x_i^{\mathrm{exper}})+e_i,
\] 
where $e_i$'s are independent $\mathrm{N}(0,\sigma^2)$ noises. The $z_i^{\mathrm{female}}$ are $\pm1$-valued, where $z_i^{\mathrm{female}}=1$ indicates that the $i$th observation is a female, and $-1$ otherwise. We set $z_i^{\mathrm{married}}=1$ if the $i$th observation is married, and $-1$ otherwise.  We centralize $z_i^{\mathrm{educ}}$ and $z_i^{\mathrm{tenure}}$ before applying the partial linear model, {i.e.}, $\sum_{i=1}^nz_i^{\mathrm{educ}}=\sum_{i=1}^nz_i^{\mathrm{tenure}}=0$. The $x_i^{\mathrm{exper}}$'s are re-scaled so that they lie in $(0,1)$. To evaluate the performance of the proposed method, we use $300$ observations as the training data, with the rest of the $226$ observations left as the testing data to compute the prediction mean-squared error. The prior specification and hyperparameters for the MCMC sampler are set as follows: $\llh= 1.2$, $\uh = 2$, $m = 3$, $\pi(\sigma^2)\propto [\sigma^2]^{-2}\exp(-1/\sigma^2)$ ({i.e.}, the inverse-Gamma density), $\pi_\bbeta=\mathrm{N}(\zero,10^2\eye_4)$, and $\pi_\xi=\mathrm{N}(0,10^2)\mathbf{1}(|\xi|\leq 100)$. The range of $K$ is set to be $\{11,12,\ldots,20\}$. 

\begin{figure}[!ht] 
\centerline{\includegraphics[width=.9\textwidth]{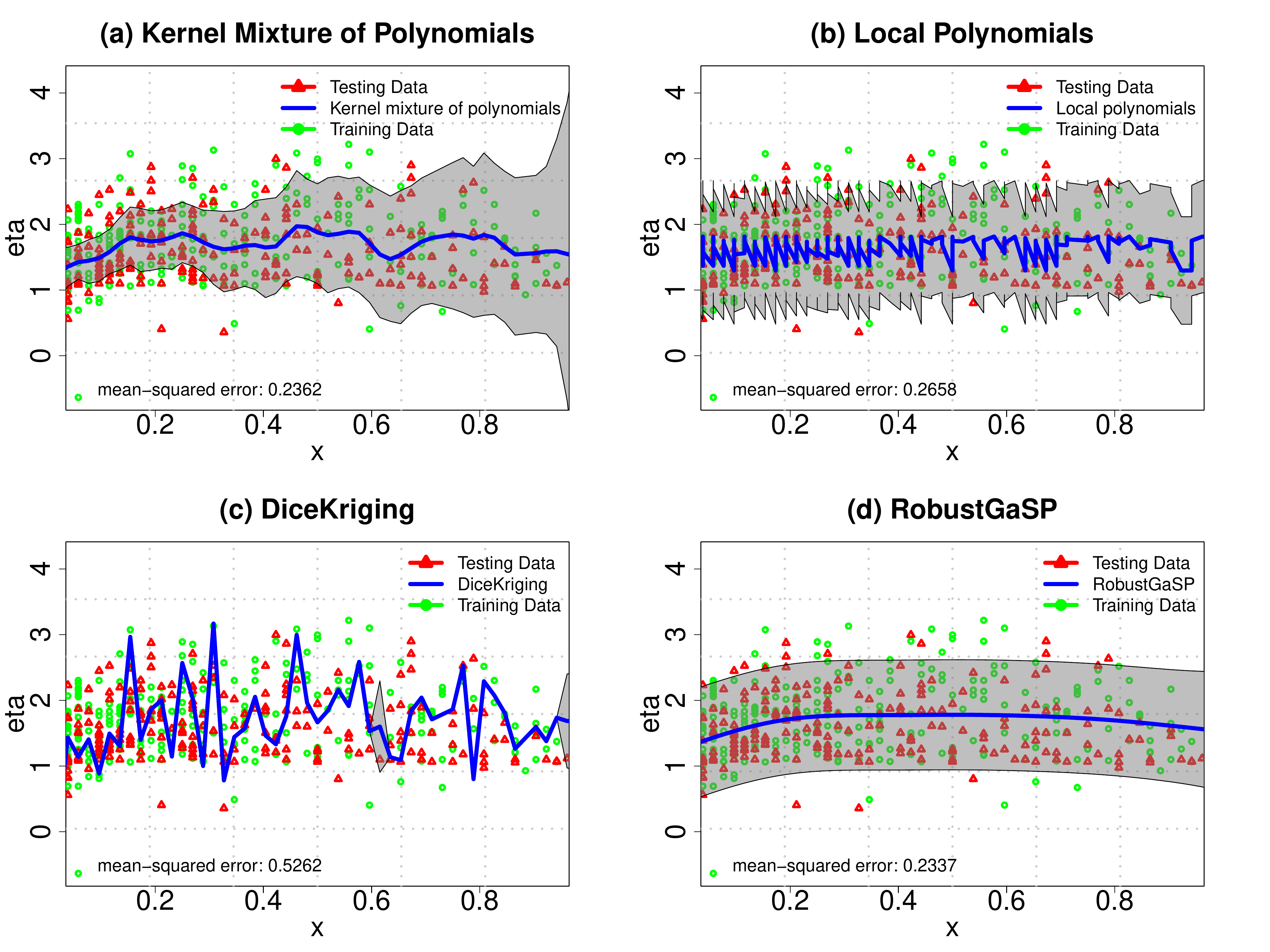}}
\caption[]{Wage data example. Shaded regions are point-wise $95\%$-credible/confidence intervals. The solid lines are point-wise posterior means/point estimators of $\eta$. Circle scatter points are the training observations, and triangle scatter points are testing responses. 
Panels (b), (c), and (d) are computed using the discrepancy data $(x_i,y_i-\bz_i\transpose\widehat{\bbeta}_{np})_{i=1}^n$, where $\widehat{\bbeta}_{np}$ is the estimate of $\bbeta$ computed using the np package. }
\label{Fig:wage_data_nonparametric_component}
\end{figure}

We calculate the posterior means and the posterior $95\%$-credible intervals for $\bbeta$. For comparison, we also provide the least-squared estimate of $\bbeta$ and the estimate 
computed by the \verb|np| package \citep{hayfield2008nonparametric}. The results 
are summarized in Table \ref{tab:real_data} (the ``KMP'' column represents the posterior means of $\bbeta$), showing that the kernel mixture of polynomials estimate is closer to the \verb|np| package estimate than the least-squared estimate, and all $3$ point estimates of $\bbeta$ lie in the posterior $95\%$-credible intervals. 
\begin{table}[!ht]
  \centering
  \caption{Wage data example: Inference of $\bbeta$}
    \begin{tabular}{ccccc}
    \hline\hline
          & KMP & $95\%$ credible intervals & np package & $\widehat{\bbeta}_{LS}$ \\
    \hline
    female & -0.1214 & (-0.2534, -0.0071) & -0.1287 & -0.0921 \\
    married & 0.0249 & (-0.0943, 0.1629) & 0.0279 & 0.3209 \\
    educ & 0.0903 & (0.0373, 0.1404) & 0.0891 & 0.1257 \\
    tenure & 0.0175 & (-0.0032, 0.0372) & 0.0167 & 0.0152 \\
    \hline\hline
    \end{tabular}%
  \label{tab:real_data}%
\end{table}%
For the nonparametric component, we compute the kernel mixture of polynomials prediction on the testing dataset. The comparison with the true testing responses is demonstrated in Panel (a) of Figure \ref{Fig:wage_data_nonparametric_component}. The $3$ alternatives based on the \verb|np| estimate of $\bbeta$ for estimating the nonparametric component $\eta$ are: the local polynomial regression, DiceKriging, and RobustGaSP. The performance of these $3$ competitors are visualized in Panels (b), (c), and (d) in Figure \ref{Fig:wage_data_nonparametric_component}, respectively. The local polynomial regression estimate does not outperform the kernel mixture of polynomials in terms of the prediction mean-squared error, and the prediction curve is highly non-smooth. DiceKriging does not work in this scenario: the prediction mean-squared error is large, the prediction curve is highly non-smooth, and the point-wise confidence intervals show singularity in estimating the covariance matrix. RobustGaSP, though gives similar prediction mean-squared error compared to the kernel mixture of polynomials, does not capture the local nonlinearity of the nonparametric component. In addition, the point-wise $95\%$ confidence/credible intervals for the local polynomial regression and RobustGaSP are wider than those given by the kernel mixture of polynomials when $x\in(0,0.6)$. Since the data are dense in the region $(0,0.6)$, the point-wise credible intervals estimated by the kernel mixture of polynomials are thinner in this region; In contrast, the design points are sparser in the region $x\geq0.6$, and correspondingly there exists larger uncertainty in estimating $\eta$. 
Namely, the uncertainty of the kernel mixture of polynomials is adaptive to the distribution of the design points.

\section{Discussion} 
\label{sec:discussion} 
In this work we assume that an upper bound $B$ for the supremum norm of the derivatives of the true regression function $f_0$ is known and is used to construct the prior for the coefficients $(\xi_{\bk\bs}:\bk\in[K]^p,|\bs|=0,\cdots,m)$. In addition, the prior for the standard deviation $\sigma$ of the noises is supported on a compact interval $[\lsigma,\usigma]$ containing $\sigma_0$. These two restrictions can be potentially inconvenient, as $B$, $\lsigma$, and $\usigma$ are typically unknown in practice. Here we provide the following non-adaptive alternative solution. The prior specification is slightly modified as follows: First $K$ is fixed at $K_n = \lceil (n/\log n)^{1/(2\alpha+p)}\rceil$, $h$ is fixed at $2/K$, and $\bmu_\bk$ is fixed at $\bmu_\bk^\star$ for all $\bk\in[K_n]^p$. Namely, we consider the following simplified kernel mixture of polynomials:
\begin{align}
\label{eqn:kernel_mix_poly_fixed_design}
f(\bx)=\sum_{\bk\in[K_n]^p}\sum_{\bs:|\bs|=0}^m\xi_{\bk\bs}\psi_{\bk\bs}(\bx),\quad\text{with}\quad
\psi_{\bk\bs}(\bx)=
\frac{\varphi_h(\bx-\bmu_\bk^\star)(\bx-\bmu_\bk^\star)^{\bs}}{\sum_{\bl\in[K_n]^p}\varphi_h(\bx-\bmu_\bl^\star)}.
\end{align}
Then the prior $\Pi^n$ on $(f,\sigma)$ is imposed as follows: $\xi_{\bk\bs}\mid\sigma\sim\mathrm{N}(0,n^2\sigma^2)$ independently for all $\bk\in[K_n]^p$ and $0\leq|\bs|\leq m$ given $\sigma$. Note here we relax the assumption that $\xi_{\bk\bs}$ is 
 supported on $[-B,B]$. 
Denote the corresponding induced conditional prior on $f$ given $\sigma$ by $\Pi_f^n(\cdot\mid\sigma)$. We impose $\sigma^2$ with the inverse-Gamma prior $\Pi_\sigma=\mathrm{IG}(a_\sigma,b_\sigma)$, with $\pi_{\sigma^2}(\sigma^2)\propto (\sigma^2)^{-a_\sigma/2-1}\exp[-b_\sigma/(2\sigma^2)]$ for some $a_\sigma,b_\sigma\geq2$. In particular, we impose the prior on $\xi_{\bk\bs}$ with a scaling factor $n$ so that for large $n$ the prior is weakly informative. 
We provide the rate of contraction under the fixed-design regression setting with respect to the empirical $L_2$-distance
$\|f-f_0\|_{L_2(\mathbb{P}_n)}=\left[(1/n)\sum_{i=1}^n(f(\bx_i)-f_0(\bx_i))^2\right]^{1/2}$
in the following theorem, the proof of which is deferred to the supplementary material. 
\begin{theorem}[Rate of contraction, fixed-design regression]\label{thm:contraction_rate_fixed_design}
  Under the setup and the prior $\Pi^n(\cdot,\cdot)=\Pi_f^n(\cdot\mid\sigma)\times\Pi_\sigma(\cdot)$ above, it holds that
  \[
  \Pi^n\left(\left\|f-f_0\right\|_{L_2(\mathbb{P}_n)}^2>M\eps_n^2\mid\calD_n\right)\to 0
  \]
  in $\mathbb{P}_0$-probability for some constant $M>0$, where $\calD_n$ is the observed data $(y_i)_{i = 1}^n$, and $\eps_n = n^{-{\alpha}/{(2\alpha+p)}}(\log n)^{\alpha/(2\alpha+p)}$. 
\end{theorem}
In general, in developing contraction rates for nonparametric regression with respect to the $\|\cdot\|_{L_2(\mathbb{P}_\bx)}$ distance, it requires significant work to drop the boundedness assumption on the space of regression functions. The underlying reason is that when the space of functions are not uniformly bounded, existence of certain test functions satisfying the so-called Ghosal \& van der Vaart conditions originally presented in \cite{ghosal2000convergence} may not exist. To the best of our knowledge, \cite{xie2017theoretical} is the only literature addressing the posterior contraction of Bayesian nonparametric regression with respect to $\|\cdot\|_{L_2(\mathbb{P}_\bx)}$ in a systematic framework without imposing uniform boundedness on function spaces. It will be interesting to extend the techniques there to study the posterior contraction of the kernel mixture of polynomials without requiring the boundedness assumption.


Another feasible extension is variable selection in high dimensions. To be more specific, suppose that the design space $\calX = [0, 1]^p$ is with high dimensionality in the sense that $p \gg n$, and the true regression function $f_0$ only depends on a subset of the covariates $\{x_{j_1},\ldots, x_{j_q}\}\subset\{x_1,\ldots, x_p\}$, where $q < n$ is the intrinsic dimension. 
To tackle such a high-dimensional variable selection problem in nonparametric regression, we modify the kernel mixture of polynomials as follows: 
Let the kernel function $\varphi(x_1,\cdots,x_p):\mathbb{R}^p\to[0,1]$ be of the product form
  $\varphi(\bx) = \varphi(x_1,\cdots,x_p) = \prod_{j=1}^p\varphi^1(x_j)$,
  where $\varphi^1:\mathbb{R}\to[0,1]$ is a univariate kernel function. We introduce the auxiliary binary variables $(z_j)_{j=1}^p$ to indicate whether the $j$th covariate is active or not, and modify the kernel mixture of polynomial system as follows:
  \[
  \psi_{\bk\bs}(x_1,\cdots,x_p\mid z_1,\cdots,z_p) = \frac{\prod_{j=1}^p[\varphi^1_h(x_j-\mu_{\bk j})]^{z_j}}{\sum_{\bl\in[K]^p}\prod_{j=1}^p[\varphi^1_h(x_j-\mu_{\bl j})]^{z_j}}\prod_{j=1}^p(x_j-\mu_{\bk j}^\star)^{s_jz_j}.
  \]
  Then by letting
  \[
  f(x_1,\cdots,x_p\mid z_1,\cdots,z_p)=\sum_{\bk\in[K]^p}\sum_{\bs:|\bs|=0}^m\xi_{\bk\bs}\psi_{\bk\bs}(x_1,\cdots,x_p\mid z_1,\cdots,z_p),
  \]
  we can see that if $z_j = 0$,  $f$ does not depend on $x_j$. In fact, the basis function $\psi_{\bk\bs}$ only depends on covariates that are active, \emph{i.e.}, $(x_j:z_j = 1,j=1,\cdots,p)$. One can further impose independent $\mathrm{Bernoulli}(p)$ prior on $(z_1,\cdots,z_p)$ and complete the hierarchical Bayesian model by adopting the prior specification in Subsection \ref{sub:prior_specification} for the rest of the parameters.
  We believe that it will be interesting to investigate the convergence properties of such a modification of the kernel mixture of polynomials.    

There are also several other potential extensions of the current work. 
Firstly, we develop the theoretical results under the assumption that the noises $(e_i)_{i=1}^n$ are Gaussian. In cases where the noises are only assumed to be sub-Gaussian, further exploration of the convergence properties can be investigated. Secondly, the design points are assumed to be random in the present paper. In cases where the design points are fixed, which is also a common phenomenon in many physical experiments \citep{tuo2015efficient}, theoretical results for the kernel mixture of polynomials can be further extended using the techniques developed for non-independent nor-identically distributed observations by \cite{ghosal2007convergence} or fixed-design nonparametric regression by \cite{xie2017theoretical}. 
Secondly, we have only considered the case where the true regression function $f_0$ is in the $\alpha$-H\"older function class. It is also interesting to extend the current framework to the case where $f_0$ is in the $\alpha$-Sobolev function space. Roughly speaking, when $\alpha$ is an integer, a function $f$ is called $\alpha$-Sobolev if the corresponding $(\alpha-1)$th derivatives are squared integrable. The almost-sure existence of derivatives of $\alpha$-Sobolev functions guarantees that the approximation lemma (Lemma \ref{lemma:prior_thickness_lemma}) is still applicable, and hence it is reasonable to expect that our theory also applies to an $\alpha$-Sobolev $f_0$. 
In addition, when applying the kernel mixture of polynomials to the partial linear model, we only consider the case where $\mathbb{E}\bz=0$ and ${\bx}$ is independent of $\bz$, indicating that the linear component and the nonparametric component are orthogonal. On one hand, the idea of orthogonality has  been explored in the literature of calibration of inexact computer models \citep{plumlee2016orthogonal,plumlee2017bayesian}, and therefore exploring the application of the kernel mixture of polynomials to calibration of orthogonal computer models is a promising extension. 
On the other hand, 
it is also interesting to investigate the convergence theory when the two components are not orthogonal to each other. 
We also expect that the kernel mixture of polynomials can be applied to other semiparametric models besides the partial linear model, such as the single-index model \citep{ichimura1993semiparametric}, the projection pursuit regression \citep{friedman1981projection}, etc. 
Finally, we have developed a theoretical support for a sieve maximum likelihood estimator with compact restrictions on the parameter spaces. In particular, the loss function is of the least-squared form. From the computational perspective, an efficient optimization technique that is also scalable to big data can be designed to obtain the frequentist estimator in light of the rich literature of solving nonlinear least-squared problems \citep{nocedal2006numerical}. 

\section*{Appendix: Notations}

  For $1\leq r\leq\infty$, we use $\|\cdot\|_r$ to denote both the $\ell_r$-norm on any finite dimensional Euclidean space and the $L_r$-norm of a measurable function. We follow  the convention that when $r=2$, the subscript is omitted, {i.e.}, $\|\cdot\|_2=\|\cdot\|$. \xx For any integer $n$, denote $[n]=\{1,\ldots,n\}$. For $x^\star\in\mathbb{R}^p$ and $\epsilon>0$, denote $B_r(\bx^\star,\epsilon)=\{\bx\in\mathbb{R}^p:\|\bx-\bx^\star\|_r<\epsilon\}$ for $1\leq r\leq \infty$. We use $\lfloor x\rfloor$ to denote the maximal integer no greater than $x$, and $\lceil x\rceil$ to denote the minimum integer no less than $x$. The notations $a\lesssim b$ and $a\gtrsim b$ denote the inequalities up to a positive multiplicative constant. 
  
  We slightly abuse the notation and do not distinguish between a random variable and its realization. Given a distribution $\mathbb{P}_\bx$ on $\calX$, we denote the $L_r(\mathbb{P}_\bx)$-norm of a measurable function $f$ by
  $\|f\|_{L_r(\mathbb{P}_\bx)}=\{\int_\calX|f(\bx)|^r\mathbb{P}_\bx(\mathrm{d}\bx)\}^{{1}/{r}}\nonumber$
  for any $r\in[1,\infty)$. The notation $\mathbf{1}(A)$ denotes the indicator of the event $A$. 
  We refer to $\calP$ as a statistical model if it consists of a class of densities on a sample space $\calX$ with respect to some underlying $\sigma$-finite measure. 
  Given a statistical model $\calP$ and the independent and identically distributed data $\calD_n=(\bx_i)_{i=1}^n$ from some $P\in\calP$, the prior and the posterior distribution on $\calP$ are always denoted by $\Pi(\cdot)$ and $\Pi(\cdot\mid\calD_n)$, respectively. 
  We use $p_\bx(\bx)$ or $p(\bx)$ to denote the density of $x$, $\mathbb{P}_\bx$ to denote the distribution of $\bx$, and $\mathbb{E}_\bx$ for the corresponding expected value. In particular, $\phi$ denotes the probability density function of the (univariate) standard normal distribution, and we use the shorthand notation $\phi_\sigma(y)=\phi(y/\sigma)/\sigma$. The Hellinger distance between two densities $p_1,p_2$ is denoted by $H(p_1,p_2)$, and the Kullback-Leibler divergence is denoted by $D_{\mathrm{KL}}(p_1||p_2)=\int p_1(\bx)\log({p_1}(\bx)/{p_2}(\bx))\mathrm{d}\bx$. 


\begin{supplement}
\sname{Supplement A}\label{suppA} 
\stitle{Supplementary Material for ``Adaptive Bayesian nonparametric regression using kernel mixture of polynomials with application to partial linear model''}
\slink[url]{}
\sdescription{The supplementary material contains proofs for Section \ref{sec:theoretical_properties_of_the_bayesian_kernel_regression} and Section \ref{sec:bayesian_kernel_smoothing_in_partial_linear_models}, and cited theorems and results. }
\end{supplement}

\bibliographystyle{ba}
\bibliography{reference}

\begin{thebibliography}{66}
\newcommand{\enquote}[1]{``#1''}
\expandafter\ifx\csname natexlab\endcsname\relax\def\natexlab#1{#1}\fi
\expandafter\ifx\csname url\endcsname\relax
  \def\url#1{{\tt #1}}\fi
\expandafter\ifx\csname urlprefix\endcsname\relax\def\urlprefix{URL }\fi
\ifx\endbibitem\undefined \let\endbibitem\relax\fi

\bibitem[{Affandi et~al.(2013)Affandi, Fox, and
  Taskar}]{affandi2013approximate}
Affandi, R.~H., Fox, E., and Taskar, B. (2013).
\newblock \enquote{Approximate inference in continuous determinantal
  processes.}
\newblock In {\em Advances in Neural Information Processing Systems\/},
  1430--1438.
\endbibitem

\bibitem[{Bhattacharya et~al.(2014)Bhattacharya, Pati, and
  Dunson}]{bhattacharya2014anisotropic}
Bhattacharya, A., Pati, D., and Dunson, D. (2014).
\newblock \enquote{Anisotropic function estimation using multi-bandwidth
  {G}aussian processes.}
\newblock {\em Annals of Statistics\/}, 42(1): 352.
\endbibitem

\bibitem[{Bhattacharya et~al.(2015)Bhattacharya, Pati, Pillai, and
  Dunson}]{bhattacharya2015dirichlet}
Bhattacharya, A., Pati, D., Pillai, N.~S., and Dunson, D.~B. (2015).
\newblock \enquote{Dirichlet-{L}aplace priors for optimal shrinkage.}
\newblock {\em Journal of the American Statistical Association\/}, 110(512):
  1479--1490.
\endbibitem

\bibitem[{Bickel et~al.(2012)Bickel, Kleijn et~al.}]{bickel2012semiparametric}
Bickel, P., Kleijn, B., et~al. (2012).
\newblock \enquote{The semiparametric {B}ernstein-von {M}ises theorem.}
\newblock {\em The Annals of Statistics\/}, 40(1): 206--237.
\endbibitem

\bibitem[{Bickel et~al.(1998)Bickel, Klaassen, Ritov, Wellner
  et~al.}]{bickel1998efficient}
Bickel, P.~J., Klaassen, C.~A., Ritov, Y., Wellner, J.~A., et~al. (1998).
\newblock \enquote{Efficient and adaptive estimation for semiparametric
  models.}
\endbibitem

\bibitem[{Cabrera(2012)}]{cabrera2012locpol}
Cabrera, J. (2012).
\newblock \enquote{locpol: {K}ernel local polynomial regression. {R} package
  version 0.4-0.}
\endbibitem

\bibitem[{Castillo and van~der Vaart(2012)}]{castillo2012needles}
Castillo, I. and van~der Vaart, A. (2012).
\newblock \enquote{Needles and straw in a haystack: {P}osterior concentration
  for possibly sparse sequences.}
\newblock {\em The Annals of Statistics\/}, 40(4): 2069--2101.
\endbibitem

\bibitem[{Celeux et~al.(2000)Celeux, Hurn, and
  Robert}]{celeux2000computational}
Celeux, G., Hurn, M., and Robert, C.~P. (2000).
\newblock \enquote{Computational and inferential difficulties with mixture
  posterior distributions.}
\newblock {\em Journal of the American Statistical Association\/}, 95(451):
  957--970.
\endbibitem

\bibitem[{Chen et~al.(1988)}]{chen1988convergence}
Chen, H. et~al. (1988).
\newblock \enquote{Convergence rates for parametric components in a partly
  linear model.}
\newblock {\em The Annals of Statistics\/}, 16(1): 136--146.
\endbibitem

\bibitem[{Choi and Woo(2015)}]{choi2015partially}
Choi, T. and Woo, Y. (2015).
\newblock \enquote{A partially linear model using a {G}aussian process prior.}
\newblock {\em Communications in Statistics-Simulation and Computation\/},
  44(7): 1770--1786.
\endbibitem

\bibitem[{De~Jonge et~al.(2010)De~Jonge, Van~Zanten et~al.}]{de2010adaptive}
De~Jonge, R., Van~Zanten, J., et~al. (2010).
\newblock \enquote{Adaptive nonparametric {B}ayesian inference using
  location-scale mixture priors.}
\newblock {\em The Annals of Statistics\/}, 38(6): 3300--3320.
\endbibitem

\bibitem[{De~Jonge et~al.(2012)De~Jonge, Van~Zanten et~al.}]{de2012adaptive}
--- (2012).
\newblock \enquote{Adaptive estimation of multivariate functions using
  conditionally {G}aussian tensor-product spline priors.}
\newblock {\em Electronic Journal of Statistics\/}, 6: 1984--2001.
\endbibitem

\bibitem[{Devroye et~al.(2013)Devroye, Gy{\"o}rfi, and
  Lugosi}]{devroye2013probabilistic}
Devroye, L., Gy{\"o}rfi, L., and Lugosi, G. (2013).
\newblock {\em A probabilistic theory of pattern recognition\/}, volume~31.
\newblock Springer Science \& Business Media.
\endbibitem

\bibitem[{Doob(1949)}]{doob1949application}
Doob, J.~L. (1949).
\newblock \enquote{Application of the theory of martingales.}
\newblock {\em Le calcul des probabilites et ses applications\/}, 23--27.
\endbibitem

\bibitem[{Engle et~al.(1986)Engle, Granger, Rice, and
  Weiss}]{engle1986semiparametric}
Engle, R.~F., Granger, C.~W., Rice, J., and Weiss, A. (1986).
\newblock \enquote{Semiparametric estimates of the relation between weather and
  electricity sales.}
\newblock {\em Journal of the American statistical Association\/}, 81(394):
  310--320.
\endbibitem

\bibitem[{Fan and Gijbels(1996)}]{fan1996local}
Fan, J. and Gijbels, I. (1996).
\newblock {\em Local polynomial modelling and its applications: monographs on
  statistics and applied probability 66\/}, volume~66.
\newblock CRC Press.
\endbibitem

\bibitem[{Fan and Li(1999)}]{fan1999root}
Fan, Y. and Li, Q. (1999).
\newblock \enquote{Root-n-consistent estimation of partially linear time series
  models.}
\newblock {\em Journal of Nonparametric Statistics\/}, 11(1-3): 251--269.
\endbibitem

\bibitem[{Friedman and Stuetzle(1981)}]{friedman1981projection}
Friedman, J.~H. and Stuetzle, W. (1981).
\newblock \enquote{Projection pursuit regression.}
\newblock {\em Journal of the American statistical Association\/}, 76(376):
  817--823.
\endbibitem

\bibitem[{Gao and Zhou(2015)}]{gao2015rate}
Gao, C. and Zhou, H.~H. (2015).
\newblock \enquote{Rate-optimal posterior contraction for sparse {PCA}.}
\newblock {\em The Annals of Statistics\/}, 43(2): 785--818.
\endbibitem

\bibitem[{Gelman et~al.(2014)Gelman, Carlin, Stern, Dunson, Vehtari, and
  Rubin}]{gelman2014bayesian}
Gelman, A., Carlin, J.~B., Stern, H.~S., Dunson, D.~B., Vehtari, A., and Rubin,
  D.~B. (2014).
\newblock {\em Bayesian data analysis\/}, volume~2.
\newblock CRC press Boca Raton, FL.
\endbibitem

\bibitem[{Ghosal et~al.(2000)Ghosal, Ghosh, and van~der
  Vaart}]{ghosal2000convergence}
Ghosal, S., Ghosh, J.~K., and van~der Vaart, A.~W. (2000).
\newblock \enquote{Convergence rates of posterior distributions.}
\newblock {\em Annals of Statistics\/}, 28(2): 500--531.
\endbibitem

\bibitem[{Ghosal and van~der Vaart(2017)}]{ghosal2017fundamentals}
Ghosal, S. and van~der Vaart, A. (2017).
\newblock {\em Fundamentals of nonparametric Bayesian inference\/}, volume~44.
\newblock Cambridge University Press.
\endbibitem

\bibitem[{Ghosal et~al.(2007{\natexlab{a}})Ghosal, van~der Vaart
  et~al.}]{ghosal2007convergence}
Ghosal, S., van~der Vaart, A., et~al. (2007{\natexlab{a}}).
\newblock \enquote{Convergence rates of posterior distributions for non-i.i.d
  observations.}
\newblock {\em The Annals of Statistics\/}, 35(1): 192--223.
\endbibitem

\bibitem[{Ghosal et~al.(2007{\natexlab{b}})Ghosal, van~der Vaart
  et~al.}]{ghosal2007posterior}
--- (2007{\natexlab{b}}).
\newblock \enquote{Posterior convergence rates of {D}irichlet mixtures at
  smooth densities.}
\newblock {\em The Annals of Statistics\/}, 35(2): 697--723.
\endbibitem

\bibitem[{Ghosal and van~der Vaart(2001)}]{ghosal2001entropies}
Ghosal, S. and van~der Vaart, A.~W. (2001).
\newblock \enquote{{Entropies and rates of convergence for maximum likelihood
  and {B}ayes estimation for mixtures of normal densities}.}
\newblock {\em Annals of Statistics\/}, 29(5): 1233--1263.
\endbibitem

\bibitem[{Gu et~al.(2016)Gu, Palomo, and Berger}]{Gu2016RobustGaSP}
Gu, M., Palomo, J., and Berger, J.~O. (2016).
\newblock \enquote{Robust {G}a{S}{P}: an {R} Package for objective {B}ayesian
  emulation of complex computer model codes.}
\newblock {\em Technical Report\/}.
\endbibitem

\bibitem[{Gu et~al.(2017)Gu, Wang, and Berger}]{Gu2017Robust}
Gu, M., Wang, X., and Berger, J.~O. (2017).
\newblock \enquote{Robust {G}aussian Stochastic Process Emulation.}
\newblock {\em arXiv preprint arXiv:1708.04738\/}.
\endbibitem

\bibitem[{Gy{\"o}rfi et~al.(2006)Gy{\"o}rfi, Kohler, Krzyzak, and
  Walk}]{gyorfi2006distribution}
Gy{\"o}rfi, L., Kohler, M., Krzyzak, A., and Walk, H. (2006).
\newblock {\em A distribution-free theory of nonparametric regression\/}.
\newblock Springer Science \& Business Media.
\endbibitem

\bibitem[{Hastie and Tibshirani(1990)}]{hastie1990generalized}
Hastie, T. and Tibshirani, R. (1990).
\newblock {\em Generalized additive models\/}.
\newblock Wiley Online Library.
\endbibitem

\bibitem[{Hayfield et~al.(2008)Hayfield, Racine
  et~al.}]{hayfield2008nonparametric}
Hayfield, T., Racine, J.~S., et~al. (2008).
\newblock \enquote{Nonparametric econometrics: The np package.}
\newblock {\em Journal of Statistical Software\/}, 27(5): 1--32.
\endbibitem

\bibitem[{Ichimura(1993)}]{ichimura1993semiparametric}
Ichimura, H. (1993).
\newblock \enquote{Semiparametric least squares ({SLS}) and weighted {SLS}
  estimation of single-index models.}
\newblock {\em Journal of Econometrics\/}, 58(1-2): 71--120.
\endbibitem

\bibitem[{Jasra et~al.(2005)Jasra, Holmes, and Stephens}]{jasra2005markov}
Jasra, A., Holmes, C.~C., and Stephens, D.~A. (2005).
\newblock \enquote{Markov chain Monte Carlo methods and the label switching
  problem in Bayesian mixture modeling.}
\newblock {\em Statistical Science\/}, 50--67.
\endbibitem

\bibitem[{Klein and Spady(1993)}]{klein1993efficient}
Klein, R.~W. and Spady, R.~H. (1993).
\newblock \enquote{An efficient semiparametric estimator for binary response
  models.}
\newblock {\em Econometrica: Journal of the Econometric Society\/}, 387--421.
\endbibitem

\bibitem[{Knapik et~al.(2011)Knapik, van~der Vaart, van Zanten
  et~al.}]{knapik2011bayesian}
Knapik, B.~T., van~der Vaart, A.~W., van Zanten, J.~H., et~al. (2011).
\newblock \enquote{Bayesian inverse problems with {G}aussian priors.}
\newblock {\em The Annals of Statistics\/}, 39(5): 2626--2657.
\endbibitem

\bibitem[{Kruijer et~al.(2010)Kruijer, Rousseau, van~der Vaart
  et~al.}]{kruijer2010adaptive}
Kruijer, W., Rousseau, J., van~der Vaart, A., et~al. (2010).
\newblock \enquote{{Adaptive {B}ayesian density estimation with location-scale
  mixtures}.}
\newblock {\em Electronic Journal of Statistics\/}, 4: 1225--1257.
\endbibitem

\bibitem[{Lenk(1999)}]{lenk1999bayesian}
Lenk, P.~J. (1999).
\newblock \enquote{Bayesian inference for semiparametric regression using a
  {F}ourier representation.}
\newblock {\em Journal of the Royal Statistical Society: Series B (Statistical
  Methodology)\/}, 61(4): 863--879.
\endbibitem

\bibitem[{Nadaraya(1964)}]{nadaraya1964estimating}
Nadaraya, E.~A. (1964).
\newblock \enquote{On estimating regression.}
\newblock {\em Theory of Probability \& Its Applications\/}, 9(1): 141--142.
\endbibitem

\bibitem[{Nocedal and Wright(2006)}]{nocedal2006numerical}
Nocedal, J. and Wright, S. (2006).
\newblock {\em Numerical optimization\/}.
\newblock Springer Science \& Business Media.
\endbibitem

\bibitem[{Pati et~al.(2014)Pati, Bhattacharya, Pillai, and
  Dunson}]{pati2014posterior}
Pati, D., Bhattacharya, A., Pillai, N.~S., and Dunson, D. (2014).
\newblock \enquote{Posterior contraction in sparse {B}ayesian factor models for
  massive covariance matrices.}
\newblock {\em The Annals of Statistics\/}, 42(3): 1102--1130.
\endbibitem

\bibitem[{Plumlee(2017)}]{plumlee2017bayesian}
Plumlee, M. (2017).
\newblock \enquote{Bayesian calibration of inexact computer models.}
\newblock {\em Journal of the American Statistical Association\/}, 1--12.
\endbibitem

\bibitem[{Plumlee and Joseph(2016)}]{plumlee2016orthogonal}
Plumlee, M. and Joseph, V.~R. (2016).
\newblock \enquote{Orthogonal {G}aussian process models.}
\newblock {\em arXiv preprint arXiv:1611.00203\/}.
\endbibitem

\bibitem[{Rasmussen and Williams(2006)}]{rasmussen2006gaussian}
Rasmussen, C.~E. and Williams, C.~K. (2006).
\newblock {\em {G}aussian processes for machine learning\/}, volume~1.
\newblock MIT press Cambridge.
\endbibitem

\bibitem[{Robinson(1988)}]{robinson1988root}
Robinson, P.~M. (1988).
\newblock \enquote{Root-N-consistent semiparametric regression.}
\newblock {\em Econometrica: Journal of the Econometric Society\/}, 931--954.
\endbibitem

\bibitem[{Ro{\v{c}}kov{\'a}(2018)}]{rovckova2018bayesian}
Ro{\v{c}}kov{\'a}, V. (2018).
\newblock \enquote{Bayesian estimation of sparse signals with a continuous
  spike-and-slab prior.}
\newblock {\em The Annals of Statistics\/}, 46(1): 401--437.
\endbibitem

\bibitem[{Roustant et~al.(2012)Roustant, Ginsbourger, and
  Deville}]{roustant2012dicekriging}
Roustant, O., Ginsbourger, D., and Deville, Y. (2012).
\newblock \enquote{Dice{K}riging, {D}ice{O}ptim: {T}wo {R} packages for the
  analysis of computer experiments by kriging-based metamodelling and
  optimization.}
\newblock {\em Journal of Statistical Software\/}, 51(1): 54p.
\endbibitem

\bibitem[{Shen et~al.(2013)Shen, Tokdar, and Ghosal}]{shen2013adaptive}
Shen, W., Tokdar, S.~T., and Ghosal, S. (2013).
\newblock \enquote{{Adaptive {B}ayesian multivariate density estimation with
  {D}irichlet mixtures}.}
\newblock {\em Biometrika\/}, 100(3): 623--640.
\endbibitem

\bibitem[{Shen and Wong(1994)}]{shen1994convergence}
Shen, X. and Wong, W.~H. (1994).
\newblock \enquote{Convergence rate of sieve estimates.}
\newblock {\em The Annals of Statistics\/}, 580--615.
\endbibitem

\bibitem[{Speckman(1988)}]{speckman1988kernel}
Speckman, P. (1988).
\newblock \enquote{Kernel smoothing in partial linear models.}
\newblock {\em Journal of the Royal Statistical Society. Series B
  (Methodological)\/}, 413--436.
\endbibitem

\bibitem[{Stone(1982)}]{stone1982optimal}
Stone, C.~J. (1982).
\newblock \enquote{Optimal global rates of convergence for nonparametric
  regression.}
\newblock {\em The annals of statistics\/}, 1040--1053.
\endbibitem

\bibitem[{Szab\'o et~al.(2015)Szab\'o, van~der Vaart, and van
  Zanten}]{szabó2015}
Szab\'o, B., van~der Vaart, A.~W., and van Zanten, J.~H. (2015).
\newblock \enquote{Frequentist coverage of adaptive nonparametric Bayesian
  credible sets.}
\newblock {\em Ann. Statist.\/}, 43(4): 1391--1428.
\newline\urlprefix\url{https://doi.org/10.1214/14-AOS1270}
\endbibitem

\bibitem[{Takeda et~al.(2007)Takeda, Farsiu, and Milanfar}]{takeda2007kernel}
Takeda, H., Farsiu, S., and Milanfar, P. (2007).
\newblock \enquote{Kernel regression for image processing and reconstruction.}
\newblock {\em IEEE Transactions on image processing\/}, 16(2): 349--366.
\endbibitem

\bibitem[{Tang et~al.(2015)Tang, Sinha, Pati, Lipsitz, and
  Lipshultz}]{tang2015bayesian}
Tang, Y., Sinha, D., Pati, D., Lipsitz, S., and Lipshultz, S. (2015).
\newblock \enquote{Bayesian partial linear model for skewed longitudinal data.}
\newblock {\em Biostatistics\/}, 16(3): 441--453.
\endbibitem

\bibitem[{Tuo and Wu(2015)}]{tuo2015efficient}
Tuo, R. and Wu, C.~J. (2015).
\newblock \enquote{Efficient calibration for imperfect computer models.}
\newblock {\em The Annals of Statistics\/}, 43(6): 2331--2352.
\endbibitem

\bibitem[{van~der Vaart and van Zanten(2007)}]{van2007bayesian}
van~der Vaart, A. and van Zanten, H. (2007).
\newblock \enquote{Bayesian inference with rescaled Gaussian process priors.}
\newblock {\em Electronic Journal of Statistics\/}, 1: 433--448.
\endbibitem

\bibitem[{van~der Vaart and Zanten(2011)}]{vaart2011information}
van~der Vaart, A. and Zanten, H.~v. (2011).
\newblock \enquote{Information rates of nonparametric {G}aussian process
  methods.}
\newblock {\em Journal of Machine Learning Research\/}, 12(Jun): 2095--2119.
\endbibitem

\bibitem[{van~der Vaart and van Zanten(2008)}]{van2008rates}
van~der Vaart, A.~W. and van Zanten, J.~H. (2008).
\newblock \enquote{Rates of contraction of posterior distributions based on
  {G}aussian process priors.}
\newblock {\em The Annals of Statistics\/}, 1435--1463.
\endbibitem

\bibitem[{van~der Vaart and van Zanten(2009)}]{van2009adaptive}
--- (2009).
\newblock \enquote{Adaptive {B}ayesian estimation using a {G}aussian random
  field with inverse {G}amma bandwidth.}
\newblock {\em The Annals of Statistics\/}, 2655--2675.
\endbibitem

\bibitem[{Watson(1964)}]{watson1964smooth}
Watson, G.~S. (1964).
\newblock \enquote{Smooth regression analysis.}
\newblock {\em Sankhy{\=a}: The Indian Journal of Statistics, Series A\/},
  359--372.
\endbibitem

\bibitem[{Wooldridge(2015)}]{wooldridge2015introductory}
Wooldridge, J.~M. (2015).
\newblock {\em Introductory econometrics: A modern approach\/}.
\newblock Nelson Education.
\endbibitem

\bibitem[{Xie et~al.(2017)Xie, Jin, and Xu}]{xie2017theoretical}
Xie, F., Jin, W., and Xu, Y. (2017).
\newblock \enquote{A Theoretical Framework for Bayesian Nonparametric
  Regression: Orthonormal Random Series and Rates of Contraction.}
\newblock {\em arXiv preprint arXiv:1712.05731\/}.
\endbibitem

\bibitem[{Xie and Xu(2017)}]{xie2017bayesian}
Xie, F. and Xu, Y. (2017).
\newblock \enquote{Bayesian Repulsive {G}aussian Mixture Model.}
\newblock {\em arXiv preprint arXiv:1703.09061\/}.
\endbibitem

\bibitem[{Xu et~al.(2016{\natexlab{a}})Xu, Mueller, and
  Telesca}]{xu2016bayesian}
Xu, Y., Mueller, P., and Telesca, D. (2016{\natexlab{a}}).
\newblock \enquote{{Bayesian inference for latent biologic structure with
  determinantal point processes (DPP)}.}
\newblock {\em Biometrics\/}, 72(3): 955--964.
\endbibitem

\bibitem[{Xu et~al.(2016{\natexlab{b}})Xu, Xu, and Saria}]{xuyanbo2016bayesian}
Xu, Y., Xu, Y., and Saria, S. (2016{\natexlab{b}}).
\newblock \enquote{A {B}ayesian Nonparametric Approach for Estimating
  Individualized Treatment-Response Curves.}
\newblock In {\em Machine Learning for Healthcare Conference\/}, 282--300.
\endbibitem

\bibitem[{Yang et~al.(2017)Yang, Bhattacharya, and Pati}]{yang2017frequentist}
Yang, Y., Bhattacharya, A., and Pati, D. (2017).
\newblock \enquote{Frequentist coverage and sup-norm convergence rate in
  {G}aussian process regression.}
\newblock {\em arXiv preprint arXiv:1708.04753\/}.
\endbibitem

\bibitem[{Yang et~al.(2015)Yang, Cheng, and Dunson}]{yang2015semiparametric}
Yang, Y., Cheng, G., and Dunson, D.~B. (2015).
\newblock \enquote{Semiparametric {B}ernstein-von {M}ises Theorem: Second Order
  Studies.}
\newblock {\em arXiv preprint arXiv:1503.04493\/}.
\endbibitem

\bibitem[{Yoo et~al.(2016)Yoo, Ghosal et~al.}]{yoo2016supremum}
Yoo, W.~W., Ghosal, S., et~al. (2016).
\newblock \enquote{Supremum norm posterior contraction and credible sets for
  nonparametric multivariate regression.}
\newblock {\em The Annals of Statistics\/}, 44(3): 1069--1102.
\endbibitem

\end{thebibliography}


\end{document}